\documentclass[12pt]{article}

\usepackage{amsmath}
\usepackage{theorem}
\usepackage{amssymb}
\usepackage{latexsym}
\usepackage{a4wide}
\usepackage{longtable}
\usepackage{epic}

\newtheorem{thm}{Theorem}
\newtheorem{lem}[thm]{Lemma}
\newtheorem{prop}[thm]{Proposition}

{\theorembodyfont{\rmfamily} \newtheorem{alg}[thm]{Algorithm}}
{\theorembodyfont{\rmfamily} \newtheorem{exa}[thm]{Example}}
\newenvironment{rem}{\noindent{\bf Remark.}}{\newline}
\newenvironment{pf}{\noindent{\bf Proof.}}{\hbox{}\hfill $\Box$}

\newcommand{\C}{\mathbb{C}}
\newcommand{\Q}{\mathbb{Q}}
\newcommand{\Z}{\mathbb{Z}}

\newcommand{\SL}{\mathrm{\mathop{SL}}}

\newcommand{\ad}{\mathrm{\mathop{ad}}}
\newcommand{\Stab}{\mathrm{\mathop{Stab}}}
\newcommand{\rank}{\mathrm{\mathop{rank}}}

\newcommand{\mf}[1]{\mathfrak{#1}}
\newcommand{\ssl}{\mathfrak{\mathop{sl}}}
\newcommand{\g}{\mathfrak{g}}
\newcommand{\myl}{\mathfrak{l}}
\newcommand{\hh}{\mathfrak{h}}
\newcommand{\s}{\mathfrak{s}}
\newcommand{\mya}{\mathfrak{a}}
\newcommand{\rr}{\mathfrak{r}}

\begin{document}

\title{Computing representatives of nilpotent orbits of $\theta$-groups}
\author{Willem A. de Graaf\\
Dipartimento di Matematica\\
Universit\`{a} di Trento\\
Italy}
\date{}
\maketitle

\begin{abstract}
We describe two algorithms for finding representatives of the nilpotent
orbits of a $\theta$-group. The algorithms have
been implemented in the computer algebra system  {\sf GAP} (inside the
package {\sf SLA}). We comment on their performance. We apply the
algorithms to study the nilpotent orbits of $\theta$-groups, where 
$\theta$ is an N-regular automorphism of a simple Lie
algebra of exceptional type.
\end{abstract}

\section{Introduction}\label{sec:1}

Let $\g$ be a simple complex Lie algebra, and let $G$ be a connected 
algebraic group with Lie algebra $\g$. Then $G$ acts on $\g$. The orbit
structure of the action of $G$ on $\g$ has been studied in detail (we refer to 
\cite{colmcgov} for an overview). In particular, the nilpotent orbits have
been classified, using a correspondence between nilpotent $G$-orbits, and
$G$-conjugacy classes of $\ssl_2$-triples. \par
Kostant and Rallis (\cite{kostant_rallis}) considered 
decompositions of the form $\g = \g_0\oplus \g_1$, 
where $\g_i$ is the eigenspace of an involution 
$\theta$ of $\g$, corresponding to the eigenvalue $(-1)^i$. Let $G_0\subset
G$ be the connected subgroup with Lie algebra $\g_0$. Then the reductive
group $G_0$ acts on $\g_1$, and again the question arises as to what its orbits
are. In \cite{kostant_rallis}, among many other things,
it is shown that there are finitely many nilpotent 
orbits, and a correspondence between nilpotent orbits and $\ssl_2$-triples,
analogous to the one for $\g$, is established. \par
This was generalised by Vinberg in the 70's (\cite{vinberg}, \cite{vinberg2}).
He considered the decomposition of $\g$ relative to an automorphism $\theta$ 
of order $m$ (or relative to a $1$-parameter group of automorphisms if 
$m=\infty$). Here 
$$\g = \g_0 \oplus \g_1\oplus \cdots \oplus \g_{m-1},$$
where $\g_i$ is the eigenspace of $\theta$ corresponding to the eigenvalue
$\omega^i$, where $\omega$ is a primitive $m$-th root of unity. Again we get
a reductive connected algebraic group $G_0$, with Lie algebra $\g_0$, acting
on $\g_1$. The group $G_0$, together with its action on $\g_1$, is called
a $\theta$-group. In the sequel we will, by a slight abuse of language,
also call $G_0$ a $\theta$-group; the action on $\g_1$ is always understood.

\begin{exa}\label{exa:1}
Let $\g = \ssl_4(\C)$; then $G=\mathrm{SL}_4(\C)$ acts on $\g$ by conjugation.
Let $\omega\in\C$ be a primitive third root of unity. Let $\theta$ be the
automorphism of order $3$ of $\g$, given by the following matrix
$$\begin{pmatrix}
1 & \omega & \omega^2 & \omega^2 \\
\omega^2 & 1 & \omega & \omega \\
\omega & \omega^2 & 1 & 1\\
\omega & \omega^2 & 1 & 1
\end{pmatrix}.$$
Here, if on position $(i,j)$ there is $\omega^k$, then $\theta(e_{i,j}) =
\omega^k e_{i,j}$, where $e_{i,j}$ is the matrix with a $1$ on position
$(i,j)$ and zeroes elsewhere. Let $h_i = e_{i,i}-e_{i+1,i+1}$. Then we see
that $\g_0$ is spanned by $h_1$, $h_2$, $h_3$, $e_{3,4}$, $e_{4,3}$; 
which means that $\g_0 \cong \ssl_2(\C) \oplus T_2$ (where $T_2$ denotes the 
subalgebra spanned by $h_1,h_2$). Furthermore, $\g_1$ is spanned by 
$e_{1,2}$, $e_{2,3}$, $e_{2,4}$, $e_{3,1}$ and $e_{4,1}$. As $\g_0$-module
(and hence as $G_0$-module) 
$\g_1$ splits as a direct sum of two $2$-dimensional modules (spanned
respectively by $e_{2,3},e_{2,4}$ and $e_{3,1},e_{4,1}$) and a $1$-dimensional
module (spanned by $e_{1,2}$). 
\end{exa}

It is of interest to study $\theta$-groups for a number of reasons. Firstly, 
they form a class of algebraic groups for which it is not completely
hopeless to list the orbits. Secondly, many interesting representations
of algebraic groups arise as $\theta$-groups (for example $\SL(9,\C)$ acting
on $\wedge^3(\C^9)$, which is studied in \cite{elashvin}).

An orbit $G_0e$ (with $e\in \g_1$) is called nilpotent if $0$ is contained
in its closure. This happens if and only if $e$ is nilpotent as an element 
of $\g$ (that is, the adjoint map $\ad_\g(e)$ is nilpotent).
The results of Vinberg show that also here there are a finite number
of nilpotent $G_0$-orbits in $\g_1$. Secondly, there is a 
correspondence between nilpotent orbits and $\ssl_2$-triples.
Moreover, in \cite{vinberg3}, \cite{vinberg2}, 
Vinberg developed a method for obtaining
the nilpotent $G_0$-orbits in $\g_1$. It is the objective of this paper
to describe algorithmic methods, that can be implemented on computer,
for this purpose.

The conjugacy classes (in $\mathrm{Aut}(\g)$) of the finite order automorphisms
of $\g$ have been classified by Kac (\cite{kac_autom}, see also 
\cite{helgason}) in terms of so-called Kac diagrams. The connected component
of the identity of $\mathrm{Aut}(\g)$ is the group of inner automorphisms
of $\g$. It is generated by the $\exp (\ad_{\g}(e))$, for $e\in\g$ 
nilpotent. We have that
an automorphism $\theta$ is inner if and only if $\g_0$ contains a Cartan
subalgebra of $\g$ (this follows, for example,  directly from 
\cite{helgason}, Chapter X, Theorem 5.15).

As said before, the purpose of this paper is to describe algorithms for
listing the nilpotent orbits of a $\theta$-group. For this, our 
computational setup is as in \cite{gra6}. In particular, we assume that 
$\g$ is given by a multiplication table relative to a Chevalley basis.
Then all structure constants in the multiplication table are integers.
From the Kac diagram of an automorphism $\theta$ it is straightforward to
compute the matrix of $\theta$ relative to the given basis of $\g$
(this follows directly from the description of the finite order automorphisms
given in \cite{helgason}, Chapter X, Theorem 5.15). Then by linear algebra we
can construct bases of $\g_0$ and $\g_1$. In this paper we will focus on the
case where $\theta$ is an inner automorphism of $\g$. This is by far the most
common case. Secondly, in this case a Cartan subalgebra of $\g_0$ will also
be a Cartan subalgebra of $\g$. 
This implies that $\g_0$ and $\g_1$ are generated by root spaces of
$\g$. So, although the matrix of $\theta$ has coefficients in $\Q(\omega)$, 
the spaces $\g_0$, $\g_1$ are defined over $\Q$. Therefore, all computations
are completely rational, i.e., the field elements that appear all lie in $\Q$.

Several listings of nilpotent orbits of $\theta$-groups have appeared 
in the literature. We mention \cite{elashvin} ($\g$ of type $E_8$, the order
of $\theta$ equal to $3$, $\g_0$ of type $A_8$), \cite{antelash} ($\g$ of
type $E_8$, the order of $\theta$ equal to $2$, $\g_0$ of type $D_8$),
\cite{gatim} ($\g$ of type $E_8$, the order of $\theta$ equal to $5$,
$\g_0$ of type $2A_4$ and  $\g$ of type $E_7$, the order of $\theta$ equal
to $3$, $\g_0$ of type $A_2+A_5$), \cite{pervushin} ($\g$ of type $E_7$,
the order of $\theta$ equal to $4$, $\g_0$ of type $A_1+2A_3$). (We remark that
in those references all orbits are classified, so also the semisimple orbits,
and the ones of mixed type). Furthermore, in \cite{doko1} and \cite{doko2}
the nilpotent orbits are listed for the cases where $\g$ is of
exceptional type, and the order of $\theta$ is $2$. 

Littelmann (\cite{litt8}) has devised an algorithm to obtain nilpotent
orbits of $\theta$-groups. His algorithm is tailored towards the case where
the $\theta$-group is defined by a $\Z$-grading. 
We note that for $\Z$-graded algebras $\g$, a Cartan subalgebra
of $\g_0$ will also be a Cartan subalgebra of $\g$ (\cite{vinberg2}, \S 1.4).
For this reason, the algorithms in this paper work unchanged for 
$\Z$-gradings; however, we will not explicitly consider these.
The main step  of Littelmann's algorithm consists of computing sets
$w(\Delta)$, where $\Delta$ is a basis of the root system, and $w$ runs through
a set of representatives the right cosets of a Weyl subgroup $W_0$ 
of the Weyl group $W$.
For this reason, the algorithm will behave in a similar manner to
our first method (described in Section \ref{sec:3}). That is, it will
work well when the index of $W_0$ in $W$ is small. 

In \cite{popov}, Popov has developed an algorithm for computing the
Hesselink strata of the nullcone of a representation of a reductive
algebraic group. If the group in question is a $\theta$-group then
this yields an algorithm for computing the nilpotent orbits. 
Experiments with an implementation of the algorithm (due to A'Campo and Popov) 
however suggest that for the special case of $\theta$-groups the more 
specialised methods of this paper perform better. (For an example, let
$\theta$ be the N-regular inner automorphism (see Section \ref{sec:6})
of order $2$ of the Lie 
algebra of type $E_6$. Then $\g_0$ is of type $A_5+A_1$, and the module
$\g_1$ is irreducible with highest weight $(0,0,1,0,0,1)$. For this case
Popov's algorithm used 831 seconds, wheras the methods of Sections 
\ref{sec:3} and \ref{sec:4} needed respectively 5 and 93 seconds.)

This paper is organised as follows. In Section \ref{sec:2} we describe two
algorithms that we need for tasks related to Weyl subgroups of Weyl groups.
In Section \ref{sec:3} we describe the first method for listing nilpotent
orbits of $\theta$-groups. It uses the classification of the nilpotent orbits 
of $\g$. For each such orbit it is decided whether it intersects with
$\g_1$, and if this is the case, the $G_0$-orbits of the intersection are
determined. This method works well if the index of the Weyl group of $\g_0$
in the one of $\g$ is small; and the method ceases to work well if this index
is large. In that case, the dimensions of the spaces $\g_0$ and $\g_1$
will be small (compared to the dimension of $\g$). 
In Section \ref{sec:4} we describe a second method, which is based
on Vinberg's theory of carrier algebras. This method works well if the spaces
$\g_0$ and $\g_1$ are small-dimensional. In this sense the two algorithms
complement each other. All algorithms have been implemented
in the computer algebra system {\sf GAP}4 (\cite{gap4}), inside the package
{\sf SLA} (\cite{sla}). In Section \ref{sec:5} we briefly report on practical 
experiences\footnote{All timings reported in this paper have been obtained
on a 2GHz machine, with 1GB of memory for {\sf GAP}.}
with these implementations. 
In Section \ref{sec:6}, we apply the algorithms to study the nilpotent 
orbits in case $\theta$ is an 
N-regular automorphism of $\g$, where $\g$ is of exceptional type, 
of orders between $2$ and the Coxeter number. 

{\bf Acknowledgements:} I would like to thank Alexander Elashvili and
Oksana Yakimova for inspiring conversations on the subject of this paper.

\section{Algorithms for Weyl subgroups}\label{sec:2}

Let Let $\Phi$ be a root system, with basis of simple roots 
$\Delta = \{\alpha_1,\ldots, \alpha_l\}$. Then the group generated
by the reflections $s_i = s_{\alpha_i}$ is the Weyl group of $\Phi$,
which we denote by $W$. Let $\beta_1,\ldots,\beta_s$ be roots in $\Phi$ 
that form a basis of a root subsystem $\Psi$ of $\Phi$. Then the 
Weyl group $W_0$ of $\Psi$ is a subgroup of $W$, generated by the reflections
$s_{\beta_i}$. The group $W_0$ is said to be a {\em Weyl subgroup} of $W$. 
In this section we describe algorithms for two tasks concerning Weyl subgroups
of Weyl groups: finding a set of representatives of the right cosets of 
$W_0$ in $W$, and checking whether two sets of weights are conjugate 
under $W_0$. This is the content of the first two subsections. In the third
subsection we describe an application to listing the $\pi$-systems in
$\Phi$, up to $W$-conjugacy.

It is straightforward to see that the root system $\Psi$ has a basis
consisting of positive roots in $\Phi$. Therefore, in the sequel we
assume that the $\beta_i$ are positive roots. 

Let $P$ denote the weight lattice of $\Phi$, spanned by the fundamental
weights $\lambda_1,\ldots,\lambda_l$. 
Let $(~,~)$ be a $W$-invariant inner product on the space spanned by $P$.
For two weights $\lambda,\mu$ we set $\langle \lambda,\mu^\vee\rangle =
2(\lambda,\mu)/(\mu,\mu)$.
Let $F$ denote the set of weights $\mu\in P$ 
with $\langle \mu,\alpha_i^\vee\rangle \geq 0$. This set is called
the fundamental Weyl chamber of $W$.
It is known that every $W$-orbit has a unique point in $F$. The fundamental
Weyl chamber $F_0$ of $W_0$ is the set of all $\mu\in P$ with 
$\langle \mu,\beta_i^\vee \rangle \geq 0$, for all $i$. 
Every $W_0$-orbit in $P$ contains a unique point in $F_0$.

For a $w\in W$ we denote by $\ell(w)$ the length of $w$, that is, the length
of a reduced expression for $w$ as a word in the generators $s_i$. Similarly,
for $u\in W_0$, $\ell_0(u)$ will denote the length of a reduced expression
of $u$ in the generators $s_{\beta_i}$.

\subsection{Listing right cosets}\label{sec:2.1}

It is known (see \cite{dyer}) that every right coset of a Weyl subgroup
of a Weyl group has a unique representative of shortest length. Also
there are algorithms known for finding a set of shortest coset representatives.
(For example the computer algebra system {\sc Magma} (\cite{magma}), and the
{\sf GAP}3 package {\sf Chevie} (\cite{chevie})
contain implementations of such an algorithm.)
However, I have not been able to find a reference for such an algorithm
in the literature. So for reasons of completeness, this section contains
a characterisation of these shortest representatives, that also yields
an algorithm to find them.  

\begin{lem}\label{lem:cos1}
Let $w\in W$ be such that $w^{-1}(\beta_i) >0$ for $1\leq i\leq s$. Then
$w$ is the unique element of shortest length in the coset $W_0w$.
Moreover, every coset $W_0u$ contains an element with this property. 
\end{lem}

\begin{pf}
We claim the following: let $v\in W$, and let $\beta$ be one of the 
$\beta_i$. Suppose that $v^{-1}(\beta)>0$. Then $\ell( s_{\beta}v ) >
\ell(v)$.

First we note that this is equivalent to $\ell( v^{-1} s_\beta ) > 
\ell(v^{-1})$. We use the fact that for $u\in W$ we have that $\ell(u)$
is equal to the number of positive roots $\alpha\in\Phi^+$ such that 
$u(\alpha) <0$. Let
$$ U = \{ \alpha \in\Phi^+ \mid s_\beta(\alpha) \in \Phi^+\}.$$
Then $s_\beta$ permutes $U$, hence the number of $\alpha\in U$ with
$v^{-1}(\alpha) < 0$ is equal to the number of $\alpha$ in $U$ with
$v^{-1}s_\beta(\alpha) < 0$. So $v^{-1}$ and $v^{-1}s_\beta$ receive the
same ``contribution'' from $U$ towards their lengths.
Set
$$V = \{ \alpha \in\Phi^+ \mid s_\beta(\alpha) \in \Phi^-\}.$$
Then $-s_\beta$ permutes $V$. This implies that the number of $\alpha\in V$
with $v^{-1}s_\beta(\alpha) < 0$ is equal to the number of $\alpha\in V$
with $v^{-1}(\alpha) >0$. We show that this last number is strictly bigger
than the number of $\alpha\in V$ with $v^{-1}(\alpha) < 0$. This then implies
the claim. 

Set $M_1 = \{\alpha\in V \mid (\mu,\alpha)<0\}$, and
$M_2 = \{\alpha\in V \mid (\mu,\alpha)>0\}$. Let $\alpha\in M_1$. 
From $s_\beta(\alpha) < 0$ we get that $\langle \alpha,\beta^\vee\rangle >0$.
But then
$$v^{-1} (s_\beta(\alpha)) = v^{-1}(\alpha) -\langle \alpha,\beta^\vee\rangle
 v^{-1}(\beta).$$
Now since $v^{-1}(\alpha) <0$ and $v^{-1}(\beta)>0$ we get that $v^{-1}s_\beta
(\alpha) <0$. It follows that $-s_\beta$ maps $M_1$ into $M_2$. But
$\beta \in M_2$ does not lie in $-s_\beta( M_1 )$. Hence $|M_2| > |M_1|$ and 
we are done.

Now suppose that $w^{-1}(\beta_i) >0$ for $1\leq i\leq s$. Let $u\in W_0$ be
such that $\ell( uw ) > \ell(w)$. Again let $\beta$ be one of the $\beta_i$,
with $\ell_0(s_\beta u) > \ell_0(u)$. The latter condition implies that
$u^{-1}(\beta)$ is a linear combination of the $\beta_i$ with non-negative
integer coefficients. Hence $w^{-1}u^{-1}(\beta) >0$. Therefore, by the claim
above, $\ell(s_\beta uw ) > \ell(uw)$. Now by induction on $\ell_0(u)$ we
see that $\ell(uw) > \ell(w)$ for all $u\in W_0$, $u\neq 1$.

Let $v\in W$, and $\beta$ one of the $\beta_i$. If $v^{-1}(\beta) <0$,
then by arguments similar to the ones used before, we prove that 
$\ell(s_\beta v ) < \ell(v)$ (this time $-s_\beta$ maps $M_2$ into $M_1$).
Hence, if $u \in W_0w$ does not have the property that $u(\beta_i) >0$
for $1\leq i\leq s$, then we can find $\beta_i$ with $\ell(s_{\beta_i}u )
< \ell(u)$. Continuing, we eventually find an element in the coset of 
minimal length.
\end{pf}

We say that a $w\in W$ is a shortest representative if it is the unique
representative of shortest length of the coset $W_0w$.

\begin{lem}\label{lem_short}
Let $w$ be a shortest representative. Write $w = w's_{\alpha_i}$, where 
$\ell(w') = \ell(w)-1$. Then $w'$ is a shortest representative.
\end{lem}

\begin{pf}
If not, then we can write $w' = w_1 w''$, with $w_1\in W_0$ and 
$\ell(w'') < \ell(w')$. Hence $w$ and $w''s_{\alpha_i}$ lie in the same right
$W_0$-coset. Hence $w''s_{\alpha_i}$ is not a shortest representative. 
Therefore we
can write $w''s_{\alpha_i} = w_2w'''$, with $w_2\in W_0$ and $w'''$ a shortest 
representative, $\ell( w''') < \ell(w''s_{\alpha_i})$.
But then $w = w_1w_2w'''$, and this implies $w=w'''$.
But $\ell(w''') \leq \ell(w'') < \ell(w')< \ell(w)$, which is a contradiction.
\end{pf}

Now let $R_k$ denote the set of shortest representatives of length $k$.
Lemmas \ref{lem:cos1} and \ref{lem_short} lead to the following algorithm for 
computing
$R_{k+1}$ from $R_k$. Initially we put $R_{k+1}=\emptyset$. Then for 
$1\leq i\leq l$ and $w\in R_k$ we do the following: if $\ell(ws_{\alpha_i})
> \ell(w)$ and $ws_{\alpha_i}(\beta_j) >0$ for $1\leq j\leq s$ then add 
$ws_{\alpha_i}$ to $R_{k+1}$. 

\begin{rem}
The implementation of this algorithm in {\sf GAP}4 appears to work well.
For example for $W_0$ of type $2A_4$ inside $W$ of type $E_8$ (this is obtained
by taking the set of simple roots, adding the lowest root, and deleting an
appropriate simple root), the {\sf GAP}4 implementation takes 2.1 seconds, 
whereas {\sf Chevie} and {\sc Magma} V2.14-11 need respectively 9.4 and 
89.7 seconds. In this example there are 48384 cosets.
\end{rem}

\subsection{Checking conjugacy}\label{sec:2.2}

Let $\Gamma_1 = \{\mu_1,\ldots,\mu_m\}$ and $\Gamma_2= \{\lambda_1,
\ldots,\lambda_m\}$ be two subsets of $P$. In this section we describe
how we can check efficiently whether there exists a $w\in W_0$ with
$w(\Gamma_1) = \Gamma_2$. For this we first focus on the problem of 
deciding whether there is a $w\in W_0$ with 
$w(\mu_i) = \lambda_i$ for $1\leq i\leq m$.

First we observe that it is straightforward to compute, for given $\mu\in P$,
a $w\in W_0$ and a $\lambda \in F_0$ with $w(\mu)=\lambda$. (Note that
$\lambda$ is uniquely determined by $\mu$.) Indeed, we find the smallest
index $i$ with $\langle \mu, \beta_i ^\vee \rangle < 0$. If there is
no such $i$ then $\mu\in F_0$ and we are done. Otherwise, we set 
$\mu_1 = s_{\beta_i}(\mu)$, and we continue with $\mu_1$ in place of
$\mu$. This algorithm terminates as $\mu_1 > \mu$ in the usual partial
order on $P$ (which is defined by $\nu < \eta$ if $\eta-\nu$ is a sum
of positive roots). Furthermore, by tracing the $s_{\beta_i}$ that we
applied we find $w$.

This means that the problem is easily solved if the sets $\Gamma_i$ have
only one element. Indeed, we compute $u,v\in W_0$ such that $u(\mu_1)$
and $v(\lambda_1)$ lie in $F_0$. If $u(\mu_1) = v(\lambda_1)$ then they are 
conjugate, and $w = v^{-1}u$ is such that $w(\mu_1) = \lambda_1$.
Otherwise they are not conjugate.

If the $\Gamma_i$ are larger, then, as a first step, we decide whether there is
a $w_1\in W_0$ with $w_1(\mu_1) = \lambda_1$. If there is no such $w_1$ then
our problem has no solution, and we stop. Otherwise we compute such a $w_1$.

Now we need an intermezzo on stabilisers.
Let $v\in W_0$ be such that $v(\lambda_1) = \lambda\in F_0$.
For a weight $\nu$ we consider its stabiliser $\Stab_{W_0}(\nu) = \{
u\in W_0 \mid u(\nu) = \nu \}$. We note that we have an isomorphism
$\sigma : \Stab_{W_0}(\lambda_1) \to \Stab_{W_0}(\lambda)$, by $\sigma(u) = 
vu v^{-1}$. Now as $\lambda\in F_0$ we have that $\Stab_{W_0}(\lambda)$ is 
generated by the $s_{\beta_i}$ with $\langle \lambda,\beta_i^\vee\rangle =0$
(cf. \cite{gra6}, Proposition 8.3.9 - there it is proved for the full
Weyl group, but the proof goes through also for Weyl subgroups).
Hence $\Stab_{W_0}(\lambda_1)$ is generated by all $v^{-1}s_{\beta_i}v = 
s_{v^{-1}\beta_i}$, where $\beta_i$ is such that $\langle \lambda,\beta_i^\vee
\rangle = 0$. We conclude that $\Stab_{W_0}(\lambda_1)$ is a Weyl 
subgroup of $W$; moreover, we can compute the reflections that generate it.

It is straightforward to see that $\Stab_{W_0}(\lambda_1)w_1$ is exactly the
set of elements of $W_0$ that send $\mu_1$ to $\lambda_1$. Now set $\mu_i' = 
w_1(\mu_i)$ for $i\geq 1$. Then by induction on the size of $\Gamma_i$ 
we can decide whether there exists a $w\in \Stab_{W_0}(\lambda_1)$ with
$w(\mu_i') = \lambda_i$ for $2\leq i\leq m$. Now if such a $w$ exists,
then $v=ww_1\in W_0$ has the property that $v(\mu_i) = \lambda_i$. Otherwise,
such a $v$ does not exist. 

Now we return to the more general problem, i.e., to decide whether
there exists a $w\in W_0$ with $w(\Gamma_1) = \Gamma_2$. We assume that
the $\mu_i$ and $\lambda_i$ are ordered in such a way that the matrix
$B_1 = ( (\mu_i,\mu_j))_{i,j=1}^m $ is equal to the matrix $B_2 = 
((\lambda_i,\lambda_j))_{i,j=1}^m $.
Then we compute all permutations $\tau$ of $\Gamma_1$ that leave $B_1$
invariant, i.e., such that
$$(\mu_i, \mu_j) = (\tau(\mu_i), \tau(\mu_j) )\text{ for  } 1\leq i,j\leq m.$$

Then for each such $\tau$ we check whether there is $w\in W_0$ with
$w\tau(\mu_i) = \lambda_i$. 

\begin{rem}
This algorithm works rather well in practice. First of all,
usually there are not many permutations $\tau$ that leave $B_1$ invariant.
Secondly, the basic operation of the algorithm is to compute a $\lambda\in
F_0$ conjugate to a given weight; and this can be done in few steps. 
We have used an implementation in {\sf GAP}4 of Dynkin's
algorithm for classifying so-called $\pi$-systems (cf. Section \ref{sec:4.3}), 
up to $W$-conjugacy, 
in the Lie algebra of type $E_8$. A list of $76$ $\pi$-systems was constructed,
and the algorithm for deciding conjugacy under the Weyl group of type $E_8$ was
called 3873 times. The total time used was about $71$ seconds. 
\end{rem}

\subsection{Dynkin's $\pi$-systems}\label{sec:4.3}

Let $\Gamma\subset \Phi$; then $\Gamma$ is called a $\pi$-system if
\begin{enumerate}
\item[C1)] for all $\alpha,\beta\in\Gamma$ we have $\alpha-\beta\not \in \Phi$,
\item[C2)] $\Gamma$ is linearly independent.
\end{enumerate}

We have that $\Gamma\subset \Phi$ 
is a basis of a root subsystem of $\Phi$ if and
only if it is a $\pi$-system. 

In \cite{dyn} Dynkin gave a neat algorithm to classify $\pi$-systems
of maximal rank (i.e., of rank equal to the rank of $\Phi$). This
works as follows. Let $\Gamma$ be a $\pi$-system, and $D\subset \Gamma$ a 
subset corresponding to a connected component of the Dynkin diagram of 
$\Gamma$. Then $D$ is a basis of a root subsystem of $\Phi$. To $D$ we add
the lowest root of that root subsystem. Secondly, we erase a root from
$D$, different from the one added.
This yields a $\pi$-system $\Gamma'$; which is said to be obtained
from $\Gamma$ by an {\em elementary transformation}. Dynkin showed that
all $\pi$-systems of maximal rank can be obtained (up to $W$-conjugacy)
from $\Delta$, by a series of elementary transformations.
So in order to get a list of all $\pi$-systems, up to $W$-conjugacy,
one does the following:

\begin{enumerate}
\item Find all maximal $\pi$-systems that can be obtained from $\Delta$ by
perfoming elementary transformations.
\item From this set erase $W$-conjugate copies, to obtain the set $M$.
\item Let $M'$ be the set obtained from $M$ by adding all subsets of
each element of $M$.
\item From $M'$ erase $W$-conjugate copies.  
\end{enumerate}

We remark that for checking whether two $\pi$-systems are $W$-conjugate, the
algorithm from Section \ref{sec:2.2} can be used. 

\begin{rem}
Using our implementation in {\sf GAP} of this algorithm, we 
have obtained the same tables for the root systems of exceptional
type as Dynkin (\cite{dyn}).
For the root systems of classical type this algorithm has been
applied by Lorente and Gruber (\cite{logru}). However, with our implementation
of the algorithm we obtained tables that for many root systems 
contain more $\pi$-systems.
\end{rem}

\section{Listing nilpotent orbits I}\label{sec:3}

Now we return to the set up of Section \ref{sec:1}. That is, $\theta$ is
an automorphism of $\g$ of order $m$, and $\g = \oplus_{i=0}^{m-1} \g_i$
is the corresponding $\Z/m\Z$-grading. We describe a method for finding
representatives of the nilpotent $G_0$-orbits in $\g_1$. 

\subsection{Preliminary lemmas}

We start with three lemmas which, for the case $m=2$,
have been proved in \cite{kostant_rallis}. The proofs for general $m$
are entirely similar, and therefore we omit them. A triple
$(h,e,f)$ of elements of $\g$, with $[h,e]=2e$, $[h,f]=-2f$, $[e,f]=h$,
is called an $\ssl_2$-triple.

\begin{lem}\label{lem1}
Let $e\in \g_1$ be nilpotent. Then there are $h\in \g_0$ and $f\in \g_{-1}$
such that $(h,e,f)$ is an $\ssl_2$-triple.
\end{lem}

The proof is the same as the first part of the proof of \cite{kostant_rallis}, 
Proposition 4. This lemma is also part of the content of \cite{vinberg2},
Theorem 1.
We call an $\ssl_2$-triple $(h,e,f)$ with the properties
of Lemma \ref{lem1} a {\em normal} $\ssl_2$-triple, in analogy to 
\cite{kostant_rallis}. The group $G_0$ acts on normal $\ssl_2$-triples by
$g\cdot (h,e,f) = (g\cdot h,g\cdot e, g\cdot f)$.

\begin{lem}\label{lem2}
Let $A$ be the set of nilpotent $G_0$-orbits in $\g_1$. Let $B$ be the set
of $G_0$-orbits of normal $\ssl_2$-triples. Let $\varphi : A\to B$ be
defined as follows: $\varphi(G_0\cdot e) = G_0\cdot (h,e,f)$, where $(h,e,f)$
is any normal $\ssl_2$-triple containing $e$. Then $\varphi$ is well-defined,
and bijective. 
\end{lem}

The proof is analogous to the second part of the proof of 
\cite{kostant_rallis}, Proposition 4.

\begin{lem}\label{lem3}
Let $(h,e,f)$ and $(h_1,e_1,f_1)$ be normal $\ssl_2$-triples. They are conjugate
under $G_0$ if and only if $h$ and $h_1$ are.
\end{lem}

Here the proof follows the one of \cite{kostant_rallis}, Lemma 4.

We say that an $h\in \g_0$ is {\em normal} if it lies in a normal
$\ssl_2$-triple $(h,e,f)$. By Lemmas \ref{lem1}, \ref{lem2} and \ref{lem3}
listing the nilpotent $G_0$-orbits in $\g_1$ is equivalent to listing
the $G_0$-orbits of normal $h\in \g_0$.

\subsection{Deciding normality}

We call an $h\in \g_0$ {\em admissible} if there are $e,f\in \g$ such that
$(h,e,f)$ is an $\ssl_2$-triple. Let $\hh_0$ be a fixed Cartan subalgebra 
of $\g_0$. In this section we describe an algorithm
for deciding whether a given admissible $h\in \g_0$ is normal.

\begin{prop}\label{prop:1}
Let $h\in \g_0$ be admissible and write $\g_i(k) = 
\{x\in\g_i \mid [h,x] = kx \}$. Set 
$$U = \{ u\in \g_1(2) \mid [\g_0(0),u] = \g_1(2)\}.$$
Then $U$ is dense in $\g_1(2)$. Moreover, an $e\in U$ lies in a 
normal $\ssl_2$-triple $(h,e,f)$ or there is no normal $\ssl_2$-triple
containing $h$.
\end{prop}

\begin{pf}
Consider the subalgebra
$$\mya = \bigoplus_{i\in \Z} \g_{i\bmod m}(2i).$$
Then $\mya$ is reductive. Indeed, let $\kappa$ be the Killing form of
$\g$. By standard arguments (i.e., analogous to \cite{jac}, Chapter IV, \S1.I) 
we see that
$\kappa(\g_i(k),\g_j(l))=0$ unless $j=-i$ and $l=-k$. 
So, since $\g$ is the direct sum of the $\g_i(k)$, 
the restriction of $\kappa$ to $\g_i(2i)\oplus \g_{-i}(-2i)$ is
nondegenerate. This implies that the restriction of $\kappa$ to $\mya$
is nondegenerate, and hence that $\mya$ is reductive (see \cite{bou}, \S 6, 
no. 4). Hence $\g_0(0)$ is reductive as well. Set $G_{0,h} = 
\{ g\in G_0 \mid g\cdot h = h\}$. Then $G_{0,h}$ is the subgroup of $G$
corresponding to $\g_0(0)$. Since $\mya$ is a $\Z$-graded reductive 
Lie algebra, $\g_1(2)$ has a finite number of $G_{0,h}$-orbits 
(see \cite{vinberg}, \S 2.6), and therefore a unique dense one. 
Denote this orbit by $U$. Let $u\in \g_1(2)$; then the $G_{0,h}$-orbit of
$u$ is dense if and only if it tangent space $[\g_0(0),u]$ is equal 
to $\g_1(2)$. So we get that $U = \{ u\in \g_1(2) \mid [\g_0(0),u] = \g_1(2)\}$.
Now let $e\in \g_1(2)$ lie in a normal $\ssl_2$-triple $(h,e,f)$. 
Let $K$ denote the Lie algebra spanned by $h,e,f$, i.e., $K$ is
isomorphic to $\ssl_2$. Then $\mya$
is a $K$-module. From $\ssl_2$-representation theory it follows that
$\ad e : \g_0(0) \to \g_1(2)$ is surjective. 
Consider the orbit $G_{0,h}\cdot e\subset \g_1(2)$; its tangent space is
$[\g_0(0),e]=\g_1(2)$. Therefore, $e\in U$. But then, since $U$ is a single
$G_{0,h}$-orbit, all
$u\in U$ lie in a normal $\ssl_2$-triple, $(h,u,f_u)$. So, if a given element
in $U$ does not lie in such a triple, then it follows that there is no
normal $\ssl_2$-triple containing $h$.
\end{pf}

This proposition immediately implies that the following algorithm
is correct.

\begin{alg}\label{alg:dec_normal}
Input: an admissible $h\in \g_0$.\\
Output: {\sc True} if $h$ is normal, {\sc False} otherwise.
\begin{enumerate}
\item Compute the spaces $\g_1(2)$ and $\g_{-1}(-2)$.
\item If $h\not\in [\g_1(2),\g_{-1}(-2)]$ then return {\sc False}.
\item By trying a few random elements find an $e\in \g_1(2)$ with
$[\g_0(0),e] = \g_1(2)$.
\item By solving a system of linear equations we decide whether there is
an $f\in \g_{-1}(-2)$ such that $(h,e,f)$ is a normal $\ssl_2$-triple.
If such an $f$ exists, then return {\sc True}, otherwise return {\sc False}.
\end{enumerate}
\end{alg}

\begin{rem}
In the third step we need to find random elements. This can be done
as follows. Let $\Omega$ be a finite set of integers (containing, say, 
all integers from 0 to $n$, for some $n>0$). 
Let $u_1,\ldots,u_s$ be a basis of $\g_1(2)$.
Then we choose randomly, independently, and uniformly $s$ elements $\alpha_i
\in \Omega$ and form the element $e=\sum_{i=1}^s \alpha_i u_i$. If $n$ is 
large enough then $e$ will lie the $U$ from Proposition \ref{prop:1} with
high probability; so we expect to find such an $e$ within a few steps.\par
The problem is how $n$ should be chosen. In practice, usually a small 
$n$ suffices. However, one can also proceed as follows. First we choose a 
small $n$. Then every time the random $e$ does not lie in $U$ we increase $n$.
\par
Also we remark that the correctness of the output does not depend on
the way $\Omega$ is chosen. The output is always correct; only the running
time depends on the choice of $\Omega$.
\end{rem}

\begin{rem}
It is possible to dispense with the second step. However, practical 
experience has shown that the algorithm on the average becomes more
efficient when we include it.
\end{rem}

\subsection{Finding the nilpotent orbits}\label{sec:3.3}

Now we assume that the automorphism $\theta$ is inner. Hence the Cartan
subalgebra $\hh_0$ of $\g_0$ is also a Cartan subalgebra of $\g$.
Since $\g_0$ is reductive we can write $\g_0 = \myl\oplus\rr$, where $\rr$ is
the centre of $\g_0$ and $\myl = [\g_0,\g_0]$ is semisimple. 
Note that $G_0 = L\times R$, with $L$ and $R$ the subgroups corresponding to
$\myl$ and $\rr$ respectively. The action of $R$ on $\g_0$ is trivial; so 
the $G_0$-action on $\g_0$ has the same orbits as the $L$-action.
Set $\hh_0^\myl = \hh_0 \cap \myl$; then $\hh_0^\myl$ is a Cartan
subalgebra of $\myl$. (Indeed, it is a maximal toral subalgebra of $\myl$,
hence a Cartan subalgebra, as $\myl$ is semisimple.)
Let $W_\myl = N_L(\hh_0^\myl)/C_L(\hh_0^\myl)$ be the Weyl group of $\myl$. 
This group acts on $\hh_0^\myl$, and a closed Weyl chamber is a fundamental 
domain for this action. Let $\Psi$ be the root system of $\myl$ relative to 
$\hh_0^\myl$. Let $\Pi=\{\beta_1,\ldots, \beta_m\}$ be a set of simple roots in 
$\Psi$. Then the set of all $h\in \hh_0^\myl$ with 
$\beta_i(h)\geq 0$ for $1\leq i\leq m$ is a fundamental Weyl chamber $C_\myl$.
The next lemma implies that we can assume that the candidates for 
normal $h\in \g_0$ lie in $C_{\myl}\oplus \rr$. 

\begin{lem}\label{lem4}
Let $(h_1,e_1,f_1)$ be a normal $\ssl_2$-triple. Then there is a unique 
normal $\ssl_2$-triple $(h,e,f)$ that is $G_0$-conjugate to $(h_1,e_1,f_1)$,
and such that $h=h'+u$, with $h'\in C_\myl$ and $u\in \rr$.
\end{lem}

\begin{pf}
Observe that $h_1$ lies in a Cartan subalgebra of $\g_0$, which is
$G_0$-conjugate to $\hh_0$. So by acting with an element
of $G_0$ we find a conjugate normal $\ssl_2$-triple $(h_2,e_2,f_2)$ such that
$h_2\in \hh_0$. Write $h_2 = h_3 + u$, with $h_3\in \hh_0^\myl$ and $u\in \rr$. 
Now $h_3$ has a unique $W_\myl$-conjugate $h_4$ lying in $C_\myl$. So we find a 
normal $\ssl_2$-triple $(h_4+u,e_4,f_4)$, conjugate to the original one.
Moreover, a normal $\ssl_2$-triple $(h_5+v,e_5,f_5)$ with $h_5\in 
C_\myl$, and $v\in \rr$ and $h_5\neq h_4$ cannot be $G_0$-conjugate to 
$(h_1,e_1,f_1)$, because two elements of $\hh_0^\myl$ are $L$-conjugate 
if and only if they are $W_\myl$-conjugate. 
\end{pf}

Let $\Phi$ be the root system of $\g$ with respect to $\hh_0$. Let
$W$ be the corresponding Weyl group. Then $W$ acts on $\hh_0$. It is well-known
that the action of $W$ on $\hh_0$ can be realised as follows. Let $\Delta =
\{\alpha_1,\ldots,\alpha_l\}$ be a basis of simple roots 
of $\Phi$. Let $h_{\alpha_i}$
be the unique element of $\hh_0$ with $\alpha_j(h_{\alpha_i}) = \langle 
\alpha_j,\alpha_i^\vee\rangle$. Then for the simple reflections $s_{\alpha_i}$ we
have $s_{\alpha_i}(h) = h -\alpha_i(h) h_{\alpha_i}$. Let $C\subset \hh_0$
be the fundamental Weyl chamber, i.e., $C$ consists of all $h\in \hh_0$ with
$\alpha_i(h)\geq 0$ for $1\leq i\leq l$. Since $\hh_0$ is also a Cartan
subalgebra of $\g_0$, we have that $\Psi$ is a root subsystem of $\Phi$.
This implies that $W_\myl$ is a Weyl subgroup of $W$. We note that $C_{\myl}
\oplus \rr$ is the fundamantal Weyl chamber corresponding to the action of
$W_\myl$ on $\hh_0$. We choose the positive roots in $\Psi$ such that they
are also positive in $\Phi$. Then $C\subset C_{\myl}\oplus \rr$. 
Since $W_{\myl}$ is a Weyl subgroup of $W$, by using the algorithm
of Section \ref{sec:2.1}, we can compute a set $\{w_1,\ldots,w_m\}$ of 
shortest length right coset representatives of $W_\myl$ in $W$.

Now we have the following algorithm.

\begin{alg}\label{alg:normal_list}
Input: an admissible $h\in C$.\\
Output: the list of all normal $\ssl_2$-triples $(h',e',f')$ with 
$h'\in C_\myl\oplus \rr$, and such that $h'$ is $W$-conjugate to $h$.
\begin{enumerate}
\item Set $\mathcal{H} = \{ w_i(h) \mid 1\leq i\leq m\}$. (Here the
$w_i$ are the shortest length right coset representatives, as above.)
\item Set $T=\emptyset$. 
For each $h\in \mathcal{H}$ we test whether $h$ is normal
using Algorithm \ref{alg:dec_normal}. If this is the case then
we find a normal $\ssl_2$-triple containing $h$. We add it to $T$.
\item Return $T$.
\end{enumerate}
\end{alg}

\begin{prop}\label{prop:corr}
Algorithm \ref{alg:normal_list} is correct.
\end{prop}

\begin{pf}
Note that by Lemma \ref{lem:cos1}, the right coset representatives $w_j$
are exactly the elements of $W$ that map $C$ into $C_\myl\oplus \rr$. 
So $\mathcal{H}$ is exactly the set of $h'\in C_\myl
\oplus \rr$ that are $W$-conjugate to $h$. 
From this set we take the normal elements. 
\end{pf}

\begin{rem}
Note that it is straightforward to find an $\ssl_2$-triple containing
a normal element of $\hh_0$. In fact, Algorithm \ref{alg:dec_normal} 
finds one.
\end{rem}

The nilpotent $G$-orbits in $\g$ have been classified in terms of so-called
weighted Dynkin diagrams (see \cite{cart}, \cite{colmcgov}). From such a
diagram it is straightforward to find a $h\in C$ lying in an 
$\ssl_2$-triple $(h,e,f)$ (cf. \cite{gra14}).
This also implies that $e$ is a representative of the corresponding
nilpotent $G$-orbit. Let $h_1,\ldots,h_t$ denote the 
elements of $C$ so obtained.

\begin{lem}\label{lem5}
Let $h\in \hh_0$ be normal. Then
there is a unique $h_i$ such that $h$ is $W$-conjugate to $h_i$. 
\end{lem}

\begin{pf}
The $h_i$ lie in $\ssl_2$-triples $(h_i,e_i,f_i)$. 
Let $(h,e,f)$ be a normal $\ssl_2$-triple. Then there is $e_i$ such
that $e$ is $G$-conjugate to $e_i$. But that implies that the triples
$(h,e,f)$ and $(h_i,e_i,f_i)$ are $G$-conjugate. Hence $h$ and $h_i$
are $G$-conjugate. Now two elements of $\hh_0$ are $G$-conjugate if and
only if they are $W$-conjugate. 
\end{pf}

Now we apply
Algorithm \ref{alg:normal_list} to each of the $h_i$, obtaining sets
of $T_i$. Let $T$ be their union. We claim that for each nilpotent
$G_0$-orbit in $\g_1$ the set $T$ contains exactly one $\ssl_2$-triple
$(h,e,f)$ such that $e$ is a representative of the orbit. Indeed,
let $e'\in \g_1$ be nilpotent, lying in a normal $\ssl_2$-triple
$(h',e',f')$. Then $h'$ is $G_0$-conjugate to an $h''\in \hh_0$
(as $h'$ lies in a Cartan subalgebra of $\g_0$, which is $G_0$-conjugate
to $\hh_0$). By the previous lemma $h''$ is $W$-conjugate to exactly 
one $h_i$. But by Lemma \ref{lem4}, $h''$ is $W_0$-conjugate to a unique
$h\in C_\myl\oplus \rr$. Hence Proposition \ref{prop:corr} implies that
$T_i$ contains an $\ssl_2$-triple $(h,e,f)$. We conclude that $T_i$
contains at least one $\ssl_2$-triple $(h,e,f)$ such that $h$ is 
$G_0$-conjugate to $h'$. If there were another one, $(h_1,e_1,f_1)$, then
$h$ and $h_1$ would be $G_0$-conjugate, and since they both lie in
$C_\myl\oplus \rr$, they would be equal by Lemma \ref{lem4}.

\begin{rem}
The algorithm also works when $\theta$ is an outer automorphism. However,
in that case $W_\myl$ is no longer a subgroup of $W$; therefore in Algorithm
\ref{alg:normal_list}, the set $\mathcal{H}$ needs to be the intersection
$W\cdot h\cap C_\myl\oplus \rr$. The need to compute the entire orbit $W\cdot
h$ makes the algorithm much less efficient.
\end{rem}

\begin{rem}
In \cite{doko1} also a method based on computing the set $\mathcal{H}$, as
in Algorithm \ref{alg:normal_list} is developed. However, in order to decide
whether a given $h'\in C_\myl\oplus \rr$ is normal, a more complicated procedure
is used.
\end{rem}

\section{Listing nilpotent orbits II}\label{sec:4}

In this section we describe an algorithm for finding representatives of
nilpotent orbits using Vinberg's theory of carrier algebras. Again we
divide the section into various subsections. 
We recall that $\g = \oplus_{i=0}^{m-1} \g_i$ is the grading of $\g$
relative to the automorphism $\theta$. 
Throughout we
let $\hh_0$ be a fixed Cartan subalgebra of $\g_0$.
By $W_0$ we denote 
the Weyl group of the root system of $\g_0$ (relative to $\hh_0$). 

\subsection{$\Z$-graded subalgebras}\label{sec:4.1}

We consider semisimple $\Z$-graded subalgebras $\s$ of $\g$, where
$\s_k\subset \g_{ k\bmod m}$. The following terminology is used:
\begin{itemize}
\item $\s$ is called {\em regular} if it is normalised by a Cartan subalgebra
of $\g_0$.
\item A regular $\s$ is called {\em standard} if it is normalised by $\hh_0$.
\item A regular $\s$ is called {\em complete} if it is not a $\Z$-graded
subalgebra of a $\Z$-graded regular semisimple subalgebra of the same rank.
\item $\s$ is called {\em locally flat} if $\dim \s_0 = \dim \s_1$.
\end{itemize} 

Let $\s$ be a standard $\Z$-graded semisimple subalgebra. Then $\s$ is
contained in a unique complete standard $\Z$-graded semisimple subalgebra, 
of the same rank (\cite{vinberg2}, Proposition 3). We call it the 
{\em completion} of $\s$. 

Let $e'\in \g_1$ be nilpotent; then after raplacing $e'$ by a $G_0$-conjugate
$e$ we get that $e$ lies in a complete standard locally flat $\Z$-graded
subalgebra $\s$, called the {\em carrier algebra} of $e$.
Moreover, $e\in \s_1$ is in general position, which means that
$[\s_0,e]=\s_1$ (\cite{vinberg2}, Theorem 2). 

Let $e,e'\in \g_0$ be nilpotent elements lying in carrier algebras
$\s$, $\s'$. As $\s$, $\s'$ are normalised by $\hh_0$, their respective
root systems $\Psi$, $\Psi'$ are sets of weights of $\hh_0$; therefore,
$W_0$ acts on their elements. We have that $e$, $e'$ are $G_0$-conjugate
if and only if $\Psi$, $\Psi'$ are $W_0$-conjugate (\cite{vinberg2},
Corollary to Theorem 2).

\begin{exa}\label{exa:2}
We consider the set up of Example \ref{exa:1}. Let $e=e_{2,3}+e_{3,1}$.
Let $\s_i$ for $i=-2,-1,0,1,2$ respectivey be the spaces spanned
by $\{e_{1,2}\}$, $\{e_{3,2},e_{1,3}\}$, $\{e_{1,1}-e_{3,3},e_{2,2}-e_{3,3}\}$,
$\{e_{2,3},e_{3,1}\}$, and $\{e_{2,1}\}$. Then $\s = \s_{-2}\oplus \cdots
\oplus \s_2$ is the carrier algebra of $e$. 
\end{exa}

\subsection{The completion}\label{sec:4.2}

Let $\s\subset \g$ be a standard semisimple $\Z$-graded 
subalgebra. Set $\mf{t} = \hh_0\cap \s_0$; then $\mf{t}$ is a maximal torus
in $\s_0$. Indeed, if $\mf{t}$ could be enlarged to a bigger torus, then
so could $\hh_0$, but the latter is a maximal torus in $\g_0$. It follows that
$\mf{t}$ is a Cartan subalgebra of $\s_0$, and therefore of $\s$ (as $\s$
is $\Z$-graded).  
Let $\Pi$ be a basis of simple roots of the root system of $\s$ with respect to 
$\mf{t}$. Every root has a degree:
$\deg(\alpha) = k$ if $\s_\alpha \subset \s_k$. We can choose $\Pi$ so that
$\deg(\alpha)\geq 0$ for all $\alpha \in \Pi$. Furthermore, if $\s$ is 
locally flat then
$\deg(\alpha)\in \{0,1\}$ for $\alpha\in\Pi$ (\cite{vinberg2}, \S 4.2).

Let $h_0\in \mf{t}$ be such that $\alpha(h_0) = \deg(\alpha)$ for 
$\alpha\in \Pi$. Then
$$\s_k = \{ x \in \s \mid [h_0,x] = kx \}.$$
Therefore, $h_0$ is called a {\em defining element} of $\s$. It is uniquely 
defined by the choice of $\mf{t}$ and $\Pi$. 

Set
\begin{equation}\label{eqn:1}
\mf{z} = \{ u \in \hh_0 \mid \alpha(u) = 0 \text{ for all } \alpha\in\Pi\}.
\end{equation}
Note that since $\s$ is regular, the elements of $\Pi$ are weights
of $\hh_0$, in particular we can view $\Pi$ as a subset of $\hh_0^*$.
Therefore, the definition of $\mf{z}$ makes sense. Furthermore, 
$\hh_0 = \mf{t} \oplus \mf{z}$. Indeed, $\mf{t}\cap \mf{z} = 0$,
$\dim \mf{z} + \dim \mf{t} = \dim \hh_0$. 

Now define the $\Z$-graded subalgebra $\mf{u}$ by 
$$\mf{u}_k = \{ x \in \g_{k \bmod m} \mid [\mf{z},x] = 0 \text{ and } 
[h_0,x] = kx \}.$$

\begin{prop}\label{prop:2}
$\mf{u} = \s' \oplus \mf{z}$, where $\s'$ is a standard complete $\Z$-graded
subalgebra of the same rank as $\s$, and containing $\s$. In other words,
$\s'$ is the completion of $\s$. 
\end{prop}

\begin{pf}
(cf. \cite{vinberg2}, Proposition 2).
Note that $[\mf{z},\s] = 0$, and hence $\s \subset \mf{u}$. Furthermore,
$\mf{t}\oplus \mf{z}$ is also a Cartan subalgebra of $\mf{u}_0$. Hence
the rank of $\mf{u}$ is equal to the rank of $\g_0$. By \cite{vinberg2},
Lemma 2, $\mf{u}$ is reductive. Hence $\mf{z}$ is the centre of $\mf{u}$. So
$\s' = [\mf{u}, \mf{u}]$ is semisimple, contains $\mf{s}$, and is of the same 
rank as $\s$. Moreover, by \cite{vinberg2}, Proposition 1,
$\s'$ is complete. Finally, as $\hh_0 \subset \mf{u}$, also $\s'$ is standard.
\end{pf}

Now for the remainder of this section we suppose that the automorphism
$\theta$ is inner. Then $\hh_0$ is also a Cartan subalgebra of $\g$.
Therefore, $\s$ is spanned by root spaces of $\g$.
By $\Phi$ we denote the root system of $\g$ with respect to $\hh_0$. 
Furthermore, for $\alpha\in\Phi$, $\g_\alpha$ will be the corresponding root 
space.\par
Set
\begin{align}
\Psi_0 &= \{ \alpha \in \Phi \mid \g_\alpha \subset \g_0 \text{ and }
\alpha(h_0) = 0, \alpha(\mf{z}) = 0 \},\nonumber \\
\Psi_1 &= \{ \alpha\in \Phi \mid \g_\alpha\subset \g_1 \text{ and }
\alpha(h_0) =1, \alpha(\mf{z}) = 0 \}.\label{eqn:2}
\end{align}
Let $\s'$ be the completion of $\s$.
Then from Proposition \ref{prop:2} it follows that
$\s'_0$ is the sum of $\mf{t}$ and the $\g_\alpha$ for $\alpha\in
\Psi_0$. Also, $\s_1'$ is the sum of the $\g_\alpha$ for $\alpha\in \Psi_1$.
So $\s'$ is locally flat if and only if $|\Pi|+|\Psi_0| = |\Psi_1|$.
Moreover, $h_0$ is also a defining element of $\s'$.

We summarise the above in the following algorithm. The input will be a  
basis of simple roots $\Pi\subset \Phi$ of the root system of 
a standard semisimple $\Z$-graded subalgebra $\s$. We assume that $\deg(\alpha)
\in \{0,1\}$ for $\alpha\in\Pi$ and that
we have given a decomposition $\Pi = \Pi_0\cup \Pi_1$, where $\Pi_i$ contains
the roots of degree $i$.
The output will be a defining element $h_0$ of the completion 
$\s'$ if $\s'$ is locally flat. Otherwise the output will be {\sf fail}.
The algorithm takes the following steps:

\begin{enumerate}
\item Compute a basis of $\mf{t} = \hh_0 \cap \s_0$. (This can for instance be
done by taking root vectors $x_\alpha\in\g_\alpha$, 
$x_{-\alpha}\in \g_{-\alpha}$ for $\alpha\in\Pi$;
then $\mf{t}$ will be spanned by the $[x_\alpha,x_{-\alpha}]$.)
\item Compute $h_0\in\mf{t}$, with $\alpha(h_0) = \deg(\alpha)$ for
$\alpha\in \Pi$.
\item Compute a basis of $\mf{z}\subset \hh_0$, where $\mf{z}$ is as
in (\ref{eqn:1}).
\item Compute $\Psi_0$ and $\Psi_1$, as in (\ref{eqn:2}).
\item If $|\Pi|+|\Psi_0| \neq |\Psi_1|$ then return {\sf fail}.
Otherwise return $h_0$.
\end{enumerate}

\begin{exa}
Let $\g$, $\s$ be as in Examples \ref{exa:1}, \ref{exa:2}. After some small
calculations we get $h_0=-e_{1,1}+e_{1,2}$ and $\mf{z}$ is spanned by
$e_{1,1}+e_{2,2}+e_{3,3}-3e_{4,4}$. Using this it is straightforward to see
that the completion of $\s$  is $\s$ itself. 
\end{exa}

\subsection{Obtaining the nilpotent orbits}\label{sec:4.4}

On the basis of what was said in the previous subsections
we describe an algorithm for listing the nilpotent $G_0$-orbits in $\g_1$. 
Also in this subsection we assume the automorphism $\theta$ to be inner.
Let $\Phi_0$ (respectively, $\Phi_1$) be the set of roots $\alpha$ of $\g$
such that $\g_\alpha \subset \g_0$ (respectively, $\subset \g_1$). We note that
$\Phi_0$ is a root sub-system of $\Phi$. Let
$\Delta_0$ be a basis of simple roots of $\Phi_0$. 

In the first step in our algorithm we obtain a set $\mathcal{P}$ of
$\pi$-systems (cf. Section \ref{sec:4.3})
contained in $\Phi_0\cup \Phi_1$. This set is required
to contain, up to $W_0$-conjugacy, bases of the root systems 
of all locally flat,
standard, semisimple $\Z$-graded subalgebras of $\g$. However, the set
may be bigger than necessary, i.e., it may include $W_0$-conjugate
$\pi$-systems, and it may have $\pi$-systems that are bases of root systems of
$\Z$-graded subalgebras that are not locally flat.

First we note the following. Let $\Pi = \Pi_0\cup \Pi_1$ be 
a basis of the root system of a standard semisimple $\Z$-graded
subalgebra of $\g$. Here $\Pi_i\subset \Phi_i$ contains the roots of 
degree $i$. Then $\Pi_0$ is a $\pi$-system in $\Phi_0$. This means that
the following steps will produce a suitable list $\mathcal{P}$:

\begin{enumerate}
\item Use the algorithm in Section \ref{sec:4.3} to get the set 
$\mathcal{P}_0$ of $\pi$-systems in $\Phi_0$, up to $W_0$-conjugacy.
\item For all $\Pi_0\in \mathcal{P}_0$ we find all maximal $\Pi_1\subset
\Phi_1$ such that $\Pi_0\cup\Pi_1$ is a $\pi$-system. This yields a set
$\mathcal{P}_1'$ of $\pi$-systems in $\Phi_0\cup\Phi_1$.
\item From $\mathcal{P}_1'$ we erase $W_0$-conjugate copies (see Section 
\ref{sec:2.2}), to get a set $\mathcal{P}_1''$. 
\item Finally, set $\mathcal{P}=\emptyset$. For all $\Pi_0\cup \Pi_1$
in $\mathcal{P}_1''$ we add all $\Pi_0\cup \Pi_1'$ to $\mathcal{P}$,
where $\Pi_1'$ runs through the subsets of $\Pi_1$.
\end{enumerate}

In the second step of our algorithm we find a list of normal $\ssl_2$-triples
corresponding to the nilpotent $G_0$-orbits in $\g_1$.
This step consists in running through the list 
$\mathcal{P}$, where for every $\Pi_0\cup\Pi_1$ we do the following.
First we apply the algorithm of Section \ref{sec:4.2}. If the result is 
{\sf fail}, then we go to the next element of $\mathcal{P}$. Otherwise
the result is a defining element $h_0\in\hh_0$ of the completion $\s'$
of the $\Z$-graded semisimple subalgebra of $\g$ with root system spanned
by $\Pi_0\cup \Pi_1$. Furthermore, $\s'$ is locally flat.
As outlined in Section \ref{sec:4.1}, $\s'$ corresponds
to a nilpotent $G_0$-orbit in $\g_1$. We obtain this orbit as follows.
We observe that $\tilde{h} = 2h_0$ is such that there is a normal 
$\ssl_2$-triple $(\tilde{h},\tilde{e},\tilde{f})$, where $\tilde{e}$ is 
a representative of the orbit we are after (\cite{vinberg2}, \S 4.2).
Using the notation of Section \ref{sec:3.3}, we compute the unique 
$h\in C_\myl\oplus \rr$ that is $W_0$-conjugate to $\tilde{h}$. There is 
a normal $\ssl_2$-triple $(h,e,f)$, where $e$ lies in the orbit we 
are interested in. If the element $h$ has not occurred before, then
we find such an $\ssl_2$-triple, and add it to our list. Otherwise,
we go to the next element of $\mathcal{P}$.

\begin{rem}
On the algorithms of this section some variations are possible.
Firstly, in the set $\mathcal{P}$ we can erase all $W_0$-conjugate
copies of $\pi$-systems, not just of the maximal ones. This yields a 
much smaller set, which makes the second step easier to execute.
However, practical experiences show that this leads to a less
efficient algorithm (i.e., there is much more work involved in
checking $W_0$-conjugacy for all pairs of elements of $\mathcal{P}$,
than in computing the elements $h\in C_\myl\oplus \rr$). The advantage of
only checking the maximal $\pi$-systems for $W_0$-conjugacy, is
that with each instance also a lot of subsystems are shown to be
$W_0$-conjugate.

Secondly, the algorithm of Section \ref{sec:4.2} could just return the
$h_0$, and not worry about flatness. In that case, if the corresponding
graded subalgebra is not locally flat, then $2h_0$ may not correspond to
a nilpotent orbit. This can then be checked using the methods of Section
\ref{sec:3}. However, practical experience shows that this leads to a
less efficient algorithm as well. 
\end{rem}

\begin{rem}
It is straightforward to extend the algorithm to outer autmorphisms. Indeed,
in that case, $\Phi_0$ is the same as before, wheras $\Phi_1$ becomes the
set of weights of $\hh_0$ in $\g_1$.   
\end{rem}

\section{Practical experiences}\label{sec:5}

The algorithms of this paper have been implemented in the language of 
{\sf GAP}4. In this section we report on their running times on some
sample inputs. As input
we have taken a few N-regular automorphisms of the Lie algebra of type $E_7$
(see Section \ref{sec:6}). The running times are displayed in Table 
\ref{tab:orb_class}. The set $\mathcal{H}$ is the union of all 
sets, denoted with the same letter, that occur in the first step
of Algorithm \ref{alg:normal_list} (recall that this algorithm
is called several times). This means that $\mathcal{H}$ is the set of
all admissible elements of $\hh_0$ that are tested for normality.
Secondly, the set $\mathcal{P}$ is as in Section \ref{sec:4.4}.

\begin{table}[htb]
\begin{center}
\begin{tabular}{|r|r|r|r|r|r|r|r|r|}
\hline
& \multicolumn{3}{c|}{Method I} & \multicolumn{2}{c|}{Method II} &
\multicolumn{3}{c|}{}\\
\hline
$|\theta|$ & $|W_\myl\setminus W|$ & $|\mathcal{H}|$ & time & 
$|\mathcal{P}|$ &  time & $\dim \g_0$ & $\dim \g_1$ & \# orbits \\
\hline
2 & 72 & 721 & 41 & & $\infty$ & 63 & 70 & 94\\
3 & 672 & 4627 & 83 & 4227 & 529 & 43 & 45 & 75 \\
4 & 4032 & 22939 & 475 & 4014 & 165 & 33 & 35 & 113 \\
5 & 10080 & 52109 & 1650 & 2494 & 31 & 27 & 27 & 82 \\
6 & 40320 &       & $\infty$ & 4302 & 50 & 21 & 24 & 233 \\
\hline
\end{tabular}
\end{center}
\caption{Running times of the algorithms, with input an N-regular automorphism
$\theta$ of the Lie algebra of type $E_7$. The first column has the order of
$\theta$. The next three columns display data relative to the algorithm of
Section \ref{sec:3}: the second column has the index of $W_\myl$ in $W$,
the third column has the size of $\mathcal{H}$, and the fourth column
lists the time taken. The next two columns display data relative to the
algorithm of Section \ref{sec:4}: the fifth column has the size of the set
$\mathcal{P}$, and the sixth column the total time taken. The last three 
columns list, respectively, the dimensions of $\g_0$ and $\g_1$ and the
total number of nilpotent orbits found. All running times are in seconds.
An $\infty$ indicates that the computation did not terminate within one 
hour.}
\label{tab:orb_class}
\end{table}

From the Table we see that Method I behaves well when the dimension of the
space $\g_0$ is big, because then the index of $W_\myl$ in $W$ is small,
which results in far fewer admissible elements of $\hh_0$ that need to be
checked for normality. On the contrary, Method II behaves badly in that
case, as the set $\mathcal{P}$ gets too big. However, this method gets 
quickly better when the dimensions of the spaces $\g_0$, $\g_1$ decrease.
So indeed the two methods complement each other.

\section{N-regular automorphisms}\label{sec:6}

Let $m\geq 2$ be an integer.
By results of Antonyan and Panyushev (see \cite{panyushev2}) 
there is a unique (up to conjugacy)
inner automorphism of order $m$ such that $\g_1$ contains a regular 
nilpotent element. 
For outer automorphisms a similar statement holds (see \cite{panyushev2}).
For the Lie algebra of type $E_6$ this means that if there are outer
automorphisms of order $m$, then there is exactly one (up to conjugacy)
such that $\g_1$ contains a regular nilpotent element. 

Automorphisms such that $\g_1$ contains a regular nilpotent elements are
called {\em N-regular}.
These N-regular automorphisms have a number of special properties
(see \cite{panyushev2}), but it is not immediately 
obvious which inner automorphisms are N-regular.

By using the algorithms of the previous sections we can find the
N-regular automorphisms of order $m$. Indeed, we can list representatives
of all conjugacy classes of inner
automorphisms of order $m$. For each element of the list we obtain
the nilpotent $G_0$-orbits in $\g_1$. In particular, we find representatives
of these orbits. By computing their weighted Dynkin diagrams
we check whether a regular nilpotent element occurs among them
(we refer to \cite{elasgra} 
for an algorithm to compute weighted Dynkin diagrams).

In Tables \ref{tab:prinE6}, \ref{tab:prinE7}, \ref{tab:prinE8}, 
\ref{tab:prinF4}, and \ref{tab:prinG2}
we give the Kac diagrams of the N-regular inner automorphisms of the Lie
algebras of exceptional type, of orders between $2$ and $h-1$ (where
$h$ is the Coxeter number). Furthermore, Table \ref{tab:prinE6out} contains
the Kac-diagrams of N-regular outer automorphisms of the Lie algebra
of type $E_6$.
These Kac diagrams are extended Dynkin diagrams
with labels that define the automorphism. It turns out that in our tables
only the labels 0 and 1 occur. Therefore, we give the Kac diagram by colouring
the nodes of the extended Dynkin diagram. A black node has label 1, a 
non-black node has label 0.

The contents of the tables are as follows: the first column contains the 
order of the N-regular automorphism, and the second column its Kac diagram.
The third column lists the number of nilpotent orbits in $\g_1$. Let
$\mathcal{N}$ denote the variety of all nilpotent elements in $\g_1$.
It is known (\cite{kostant_rallis}, \cite{vinberg2})
that $\mathcal{N}$ splits in irreducible components $\mathcal{N}_i$,
where $1\leq i\leq r$.
The component $\mathcal{N}_i$ is the closure of an orbit 
$G_0 e_i$, where $e_i\in \g_1$ is
nilpotent, such that the orbit $G_0 e_i$ is of maximal possible dimension.
Therefore, all components are of the same dimension, and we can compute the
$e_i$ as follows. First we list representatives of all nilpotent orbits
in $\g_1$. Let $e$ be such a representative; then the dimension of its
orbit is equal to the dimension of $[\g_0, e]$. Hence, by simple linear
algebra, we can compute the dimension of each orbit. The $e_i$ are then 
the representatives of the orbits of maximal dimension.
In particular, we get the number
of irreducible components, and their dimension. That is the content of the 
fourth and fifth column. The last column displays the rank of $\g_1$,
that is the codimension of a nilpotent orbit of maximal dimension.
(It is also the dimension of a Cartan subspace, cf. \cite{vinberg}).

When the order of $\theta$ is $2$, then all nilpotent orbits in $\g_1$ of
maximal dimension are conjugated in the group $G$ (\cite{kostant_rallis},
Theorem 6). A similar statement fails to hold for larger orders. 
Following \cite{panyushev3}, we say that an
$N$-regular automorphism for which all nilpotent $G_0$-orbits 
in $\g_1$ of maximal dimension lie in the same $G$-orbit, is
{\em very N-regular}.
We checked which $\theta$ are very N-regular
by computing the weighted Dynkin diagrams of the orbits
$G e_i$, where $G_0 e_i$ is a nilpotent $G_0$-orbit in $\g_1$ of maximal
dimension.
If $\theta$ is {\em not} very N-regular
then the number in the fourth column has an added *.
The situation where $\mathcal{N}$ has more than one
component, and $\theta$ is very N-regular, appears to
occur very rarely.

\begin{rem}
There is a also a more direct method for constructing the N-regular
inner automorphisms of order $m$ (cf. \cite{antonyan}). Let $e\in \g$ be
a regular nilpotent element lying in an $\ssl_2$-triple $(h,e,f)$. Since $e$
is even, all eigenvalues of $\ad_{\g} h$ are even integers. Let 
$$\g = \bigoplus_{k\in \Z} \g(2k) $$
be the grading of $\g$ into $\ad_{\g} h$-eiganspaces. For $0\leq i\leq m-1$
let $\g_i$ be the sum of all $\g(2k)$ such that $k\bmod m = i$. Then 
$\g = \oplus_{i=0}^{m-1} \g_i$ is the grading of $\g$ corresponding to an
$N$-regular automorphism of order $m$. It is inner as $\g_0$ contains a 
Cartan subalgebra of $\g$. Some work is still needed to obtain the Kac
diagram from the grading.  
\end{rem}

\begin{longtable}{|c|c|c|c|c|c|}
\caption{N-regular inner automorphisms of $E_6$.}\label{tab:prinE6}
\endfirsthead
\hline
\multicolumn{6}{|l|}{\small\slshape N-regular automorphisms of $E_6$.} \\
\hline 
\endhead
\hline 
\endfoot
\endlastfoot

\hline

order & Kac diagram & \# orbits & \# components & $\dim$ & $\rank$ \\
\hline

2 & 
\begin{picture}(90,50)
  \put(3,0){\circle{6}}
  \put(23,0){\circle{6}}
  \put(43,0){\circle{6}}
  \put(63,0){\circle{6}}
  \put(83,0){\circle{6}}
  \put(43,20){\circle*{6}}
  \put(43,40){\circle{6}}
  \put(6,0){\line(1,0){14}}
  \put(26,0){\line(1,0){14}}
  \put(46,0){\line(1,0){14}}
  \put(66,0){\line(1,0){14}}
  \put(43,3){\line(0,1){14}}
  \put(43,23){\line(0,1){14}}
\end{picture}
& 37 & 1 & 36 & 4\\

3 & 
\begin{picture}(90,50)
  \put(3,0){\circle{6}}
  \put(23,0){\circle{6}}
  \put(43,0){\circle*{6}}
  \put(63,0){\circle{6}}
  \put(83,0){\circle{6}}
  \put(43,20){\circle{6}}
  \put(43,40){\circle{6}}
  \put(6,0){\line(1,0){14}}
  \put(26,0){\line(1,0){14}}
  \put(46,0){\line(1,0){14}}
  \put(66,0){\line(1,0){14}}
  \put(43,3){\line(0,1){14}}
  \put(43,23){\line(0,1){14}}
\end{picture}
& 62 & 3 & 24 & 3

\\

4 & 
\begin{picture}(90,50)
  \put(3,0){\circle*{6}}
  \put(23,0){\circle{6}}
  \put(43,0){\circle*{6}}
  \put(63,0){\circle{6}}
  \put(83,0){\circle{6}}
  \put(43,20){\circle{6}}
  \put(43,40){\circle{6}}
  \put(6,0){\line(1,0){14}}
  \put(26,0){\line(1,0){14}}
  \put(46,0){\line(1,0){14}}
  \put(66,0){\line(1,0){14}}
  \put(43,3){\line(0,1){14}}
  \put(43,23){\line(0,1){14}}
\end{picture}
& 43 & 3* & 18 & 2\\

5 & 
\begin{picture}(90,50)
  \put(3,0){\circle*{6}}
  \put(23,0){\circle{6}}
  \put(43,0){\circle*{6}}
  \put(63,0){\circle{6}}
  \put(83,0){\circle{6}}
  \put(43,20){\circle{6}}
  \put(43,40){\circle*{6}}
  \put(6,0){\line(1,0){14}}
  \put(26,0){\line(1,0){14}}
  \put(46,0){\line(1,0){14}}
  \put(66,0){\line(1,0){14}}
  \put(43,3){\line(0,1){14}}
  \put(43,23){\line(0,1){14}}
\end{picture}
& 60 & 1 & 15 & 1\\

6 & 
\begin{picture}(90,50)
  \put(3,0){\circle*{6}}
  \put(23,0){\circle{6}}
  \put(43,0){\circle*{6}}
  \put(63,0){\circle{6}}
  \put(83,0){\circle*{6}}
  \put(43,20){\circle{6}}
  \put(43,40){\circle*{6}}
  \put(6,0){\line(1,0){14}}
  \put(26,0){\line(1,0){14}}
  \put(46,0){\line(1,0){14}}
  \put(66,0){\line(1,0){14}}
  \put(43,3){\line(0,1){14}}
  \put(43,23){\line(0,1){14}}
\end{picture}
& 133 & 9* & 12 & 2\\

7 & 
\begin{picture}(90,50)
  \put(3,0){\circle*{6}}
  \put(23,0){\circle*{6}}
  \put(43,0){\circle{6}}
  \put(63,0){\circle*{6}}
  \put(83,0){\circle{6}}
  \put(43,20){\circle*{6}}
  \put(43,40){\circle{6}}
  \put(6,0){\line(1,0){14}}
  \put(26,0){\line(1,0){14}}
  \put(46,0){\line(1,0){14}}
  \put(66,0){\line(1,0){14}}
  \put(43,3){\line(0,1){14}}
  \put(43,23){\line(0,1){14}}
\end{picture}
& 53 & 1 & 11 & 0

\\

8 & 
\begin{picture}(90,50)
  \put(3,0){\circle*{6}}
  \put(23,0){\circle{6}}
  \put(43,0){\circle*{6}}
  \put(63,0){\circle*{6}}
  \put(83,0){\circle*{6}}
  \put(43,20){\circle{6}}
  \put(43,40){\circle*{6}}
  \put(6,0){\line(1,0){14}}
  \put(26,0){\line(1,0){14}}
  \put(46,0){\line(1,0){14}}
  \put(66,0){\line(1,0){14}}
  \put(43,3){\line(0,1){14}}
  \put(43,23){\line(0,1){14}}
\end{picture}
& 70 & 4* & 9 & 1\\

9 & 
\begin{picture}(90,50)
  \put(3,0){\circle*{6}}
  \put(23,0){\circle*{6}}
  \put(43,0){\circle{6}}
  \put(63,0){\circle*{6}}
  \put(83,0){\circle*{6}}
  \put(43,20){\circle*{6}}
  \put(43,40){\circle*{6}}
  \put(6,0){\line(1,0){14}}
  \put(26,0){\line(1,0){14}}
  \put(46,0){\line(1,0){14}}
  \put(66,0){\line(1,0){14}}
  \put(43,3){\line(0,1){14}}
  \put(43,23){\line(0,1){14}}
\end{picture}
& 118 & 6* & 8 & 1

\\

10 & 
\begin{picture}(90,50)
  \put(3,0){\circle*{6}}
  \put(23,0){\circle*{6}}
  \put(43,0){\circle*{6}}
  \put(63,0){\circle*{6}}
  \put(83,0){\circle*{6}}
  \put(43,20){\circle{6}}
  \put(43,40){\circle*{6}}
  \put(6,0){\line(1,0){14}}
  \put(26,0){\line(1,0){14}}
  \put(46,0){\line(1,0){14}}
  \put(66,0){\line(1,0){14}}
  \put(43,3){\line(0,1){14}}
  \put(43,23){\line(0,1){14}}
\end{picture}
& 79 & 1 & 8 & 0\\

11 & 
\begin{picture}(90,50)
  \put(3,0){\circle*{6}}
  \put(23,0){\circle*{6}}
  \put(43,0){\circle*{6}}
  \put(63,0){\circle*{6}}
  \put(83,0){\circle{6}}
  \put(43,20){\circle*{6}}
  \put(43,40){\circle*{6}}
  \put(6,0){\line(1,0){14}}
  \put(26,0){\line(1,0){14}}
  \put(46,0){\line(1,0){14}}
  \put(66,0){\line(1,0){14}}
  \put(43,3){\line(0,1){14}}
  \put(43,23){\line(0,1){14}}
\end{picture}
& 63 & 1 & 7 & 0

\\

\hline
\end{longtable}

\begin{longtable}{|c|c|c|c|c|c|}
\caption{N-regular outer automorphisms of $E_6$.}\label{tab:prinE6out}
\endfirsthead
\hline
\multicolumn{6}{|l|}{\small\slshape N-regular automorphisms of $E_6$.} \\
\hline 
\endhead
\hline 
\endfoot
\endlastfoot

\hline

order & Kac diagram & \# orbits & \# components & $\dim$ & $\rank$ \\
\hline

2 & 
\begin{picture}(130,15)
  \put(10,5){\circle{6}}
  \put(40,5){\circle{6}}
  \put(70,5){\circle{6}}
  \put(100,5){\circle{6}}
  \put(130,5){\circle*{6}}
  \put(13,5){\line(1,0){24}}
  \put(43,5){\line(1,0){24}}
  \put(72,7){\line(1,0){26}}
  \put(72,3){\line(1,0){26}}
  \put(80,1){\Large $<$}
  \put(103,5){\line(1,0){24}}
\end{picture}
& 23 & 1 & 36 & 6\\

4 & 
\begin{picture}(130,15)
  \put(10,5){\circle{6}}
  \put(40,5){\circle{6}}
  \put(70,5){\circle{6}}
  \put(100,5){\circle*{6}}
  \put(130,5){\circle{6}}
  \put(13,5){\line(1,0){24}}
  \put(43,5){\line(1,0){24}}
  \put(72,7){\line(1,0){26}}
  \put(72,3){\line(1,0){26}}
  \put(80,1){\Large $<$}
  \put(103,5){\line(1,0){24}}
\end{picture}
& 20 & 1 & 18 & 2
\\

6 & 
\begin{picture}(130,15)
  \put(10,5){\circle*{6}}
  \put(40,5){\circle{6}}
  \put(70,5){\circle{6}}
  \put(100,5){\circle*{6}}
  \put(130,5){\circle{6}}
  \put(13,5){\line(1,0){24}}
  \put(43,5){\line(1,0){24}}
  \put(72,7){\line(1,0){26}}
  \put(72,3){\line(1,0){26}}
  \put(80,1){\Large $<$}
  \put(103,5){\line(1,0){24}}
\end{picture}
& 34 & 5* & 12 & 3 \\

8 & 
\begin{picture}(130,15)
  \put(10,5){\circle*{6}}
  \put(40,5){\circle{6}}
  \put(70,5){\circle{6}}
  \put(100,5){\circle*{6}}
  \put(130,5){\circle*{6}}
  \put(13,5){\line(1,0){24}}
  \put(43,5){\line(1,0){24}}
  \put(72,7){\line(1,0){26}}
  \put(72,3){\line(1,0){26}}
  \put(80,1){\Large $<$}
  \put(103,5){\line(1,0){24}}
\end{picture}
& 22 & 3* & 9 & 1\\

10 & 
\begin{picture}(130,15)
  \put(10,5){\circle*{6}}
  \put(40,5){\circle*{6}}
  \put(70,5){\circle{6}}
  \put(100,5){\circle*{6}}
  \put(130,5){\circle{6}}
  \put(13,5){\line(1,0){24}}
  \put(43,5){\line(1,0){24}}
  \put(72,7){\line(1,0){26}}
  \put(72,3){\line(1,0){26}}
  \put(80,1){\Large $<$}
  \put(103,5){\line(1,0){24}}
\end{picture}
& 25 & 2* & 8 & 1\\

12 & 
\begin{picture}(130,15)
  \put(10,5){\circle*{6}}
  \put(40,5){\circle*{6}}
  \put(70,5){\circle{6}}
  \put(100,5){\circle*{6}}
  \put(130,5){\circle*{6}}
  \put(13,5){\line(1,0){24}}
  \put(43,5){\line(1,0){24}}
  \put(72,7){\line(1,0){26}}
  \put(72,3){\line(1,0){26}}
  \put(80,1){\Large $<$}
  \put(103,5){\line(1,0){24}}
\end{picture}
& 30 & 4* & 6 & 1\\

14 & 
\begin{picture}(130,15)
  \put(10,5){\circle*{6}}
  \put(40,5){\circle{6}}
  \put(70,5){\circle*{6}}
  \put(100,5){\circle*{6}}
  \put(130,5){\circle*{6}}
  \put(13,5){\line(1,0){24}}
  \put(43,5){\line(1,0){24}}
  \put(72,7){\line(1,0){26}}
  \put(72,3){\line(1,0){26}}
  \put(80,1){\Large $<$}
  \put(103,5){\line(1,0){24}}
\end{picture}
& 19 & 1 & 6 & 0\\

16 & 
\begin{picture}(130,15)
  \put(10,5){\circle{6}}
  \put(40,5){\circle*{6}}
  \put(70,5){\circle*{6}}
  \put(100,5){\circle*{6}}
  \put(130,5){\circle*{6}}
  \put(13,5){\line(1,0){24}}
  \put(43,5){\line(1,0){24}}
  \put(72,7){\line(1,0){26}}
  \put(72,3){\line(1,0){26}}
  \put(80,1){\Large $<$}
  \put(103,5){\line(1,0){24}}
\end{picture}
& 15 & 1 & 5 & 0\\

\hline
\end{longtable}

\begin{longtable}{|c|c|c|c|c|c|}
\caption{N-regular inner automorphisms of $E_7$.}\label{tab:prinE7}
\endfirsthead
\hline
\multicolumn{6}{|l|}{\small\slshape N-regular automorphisms of $E_7$.} \\
\hline 
\endhead
\hline 
\endfoot
\endlastfoot

\hline

order & Kac diagram & \# orbits & \# components & $\dim$ & $\rank$ \\
\hline

2 &
\begin{picture}(140,40)
  \put(3,0){\circle{6}}
  \put(23,0){\circle{6}}
  \put(43,0){\circle{6}}
  \put(63,0){\circle{6}}
  \put(83,0){\circle{6}}
  \put(63,20){\circle*{6}}
  \put(6,0){\line(1,0){14}}
  \put(26,0){\line(1,0){14}}
  \put(46,0){\line(1,0){14}}
  \put(66,0){\line(1,0){14}}
  \put(63,3){\line(0,1){14}}
\put(103,0){\circle{6}}
\put(86,0){\line(1,0){14}}
\put(123,0){\circle{6}}
\put(106,0){\line(1,0){14}}
\end{picture} 
& 94 & 2 & 63 & 7 \\

3 &
\begin{picture}(140,40)
  \put(3,0){\circle{6}}
  \put(23,0){\circle{6}}
  \put(43,0){\circle*{6}}
  \put(63,0){\circle{6}}
  \put(83,0){\circle{6}}
  \put(63,20){\circle{6}}
  \put(6,0){\line(1,0){14}}
  \put(26,0){\line(1,0){14}}
  \put(46,0){\line(1,0){14}}
  \put(66,0){\line(1,0){14}}
  \put(63,3){\line(0,1){14}}
\put(103,0){\circle{6}}
\put(86,0){\line(1,0){14}}
\put(123,0){\circle{6}}
\put(106,0){\line(1,0){14}}
\end{picture}
& 75 & 1 & 42  & 3
\\

4 &
\begin{picture}(140,40)
  \put(3,0){\circle*{6}}
  \put(23,0){\circle{6}}
  \put(43,0){\circle{6}}
  \put(63,0){\circle{6}}
  \put(83,0){\circle*{6}}
  \put(63,20){\circle{6}}
  \put(6,0){\line(1,0){14}}
  \put(26,0){\line(1,0){14}}
  \put(46,0){\line(1,0){14}}
  \put(66,0){\line(1,0){14}}
  \put(63,3){\line(0,1){14}}
\put(103,0){\circle{6}}
\put(86,0){\line(1,0){14}}
\put(123,0){\circle{6}}
\put(106,0){\line(1,0){14}}
\end{picture} 
& 113 & 1 & 33 & 2\\

5 &
\begin{picture}(140,40)
  \put(3,0){\circle*{6}}
  \put(23,0){\circle{6}}
  \put(43,0){\circle{6}}
  \put(63,0){\circle*{6}}
  \put(83,0){\circle{6}}
  \put(63,20){\circle{6}}
  \put(6,0){\line(1,0){14}}
  \put(26,0){\line(1,0){14}}
  \put(46,0){\line(1,0){14}}
  \put(66,0){\line(1,0){14}}
  \put(63,3){\line(0,1){14}}
\put(103,0){\circle{6}}
\put(86,0){\line(1,0){14}}
\put(123,0){\circle{6}}
\put(106,0){\line(1,0){14}}
\end{picture} 
& 82 & 1 & 26 & 1
\\

6 &
\begin{picture}(140,40)
  \put(3,0){\circle*{6}}
  \put(23,0){\circle{6}}
  \put(43,0){\circle{6}}
  \put(63,0){\circle*{6}}
  \put(83,0){\circle{6}}
  \put(63,20){\circle{6}}
  \put(6,0){\line(1,0){14}}
  \put(26,0){\line(1,0){14}}
  \put(46,0){\line(1,0){14}}
  \put(66,0){\line(1,0){14}}
  \put(63,3){\line(0,1){14}}
\put(103,0){\circle{6}}
\put(86,0){\line(1,0){14}}
\put(123,0){\circle*{6}}
\put(106,0){\line(1,0){14}}
\end{picture} 
& 233 & 10* & 21 & 3\\

7 &
\begin{picture}(140,40)
  \put(3,0){\circle*{6}}
  \put(23,0){\circle{6}}
  \put(43,0){\circle{6}}
  \put(63,0){\circle*{6}}
  \put(83,0){\circle{6}}
  \put(63,20){\circle{6}}
  \put(6,0){\line(1,0){14}}
  \put(26,0){\line(1,0){14}}
  \put(46,0){\line(1,0){14}}
  \put(66,0){\line(1,0){14}}
  \put(63,3){\line(0,1){14}}
\put(103,0){\circle*{6}}
\put(86,0){\line(1,0){14}}
\put(123,0){\circle{6}}
\put(106,0){\line(1,0){14}}
\end{picture} 
& 112 & 3* & 18 & 1\\

8 &
\begin{picture}(140,40)
  \put(3,0){\circle*{6}}
  \put(23,0){\circle{6}}
  \put(43,0){\circle{6}}
  \put(63,0){\circle*{6}}
  \put(83,0){\circle{6}}
  \put(63,20){\circle{6}}
  \put(6,0){\line(1,0){14}}
  \put(26,0){\line(1,0){14}}
  \put(46,0){\line(1,0){14}}
  \put(66,0){\line(1,0){14}}
  \put(63,3){\line(0,1){14}}
\put(103,0){\circle*{6}}
\put(86,0){\line(1,0){14}}
\put(123,0){\circle*{6}}
\put(106,0){\line(1,0){14}}
\end{picture} 
& 163 & 2* & 17 & 1\\

9 &
\begin{picture}(140,40)
  \put(3,0){\circle*{6}}
  \put(23,0){\circle*{6}}
  \put(43,0){\circle{6}}
  \put(63,0){\circle*{6}}
  \put(83,0){\circle{6}}
  \put(63,20){\circle{6}}
  \put(6,0){\line(1,0){14}}
  \put(26,0){\line(1,0){14}}
  \put(46,0){\line(1,0){14}}
  \put(66,0){\line(1,0){14}}
  \put(63,3){\line(0,1){14}}
\put(103,0){\circle*{6}}
\put(86,0){\line(1,0){14}}
\put(123,0){\circle{6}}
\put(106,0){\line(1,0){14}}
\end{picture}
& 132 & 4* & 14 & 1
\\

10 &
\begin{picture}(140,40)
  \put(3,0){\circle*{6}}
  \put(23,0){\circle{6}}
  \put(43,0){\circle*{6}}
  \put(63,0){\circle{6}}
  \put(83,0){\circle*{6}}
  \put(63,20){\circle*{6}}
  \put(6,0){\line(1,0){14}}
  \put(26,0){\line(1,0){14}}
  \put(46,0){\line(1,0){14}}
  \put(66,0){\line(1,0){14}}
  \put(63,3){\line(0,1){14}}
\put(103,0){\circle{6}}
\put(86,0){\line(1,0){14}}
\put(123,0){\circle*{6}}
\put(106,0){\line(1,0){14}}
\end{picture} 
& 199 & 4* & 13 & 1\\

11 &
\begin{picture}(140,40)
  \put(3,0){\circle*{6}}
  \put(23,0){\circle*{6}}
  \put(43,0){\circle{6}}
  \put(63,0){\circle*{6}}
  \put(83,0){\circle{6}}
  \put(63,20){\circle*{6}}
  \put(6,0){\line(1,0){14}}
  \put(26,0){\line(1,0){14}}
  \put(46,0){\line(1,0){14}}
  \put(66,0){\line(1,0){14}}
  \put(63,3){\line(0,1){14}}
\put(103,0){\circle*{6}}
\put(86,0){\line(1,0){14}}
\put(123,0){\circle{6}}
\put(106,0){\line(1,0){14}}
\end{picture} 
& 99 & 1 & 12 & 0\\

12 &
\begin{picture}(140,40)
  \put(3,0){\circle*{6}}
  \put(23,0){\circle{6}}
  \put(43,0){\circle*{6}}
  \put(63,0){\circle{6}}
  \put(83,0){\circle*{6}}
  \put(63,20){\circle*{6}}
  \put(6,0){\line(1,0){14}}
  \put(26,0){\line(1,0){14}}
  \put(46,0){\line(1,0){14}}
  \put(66,0){\line(1,0){14}}
  \put(63,3){\line(0,1){14}}
\put(103,0){\circle*{6}}
\put(86,0){\line(1,0){14}}
\put(123,0){\circle*{6}}
\put(106,0){\line(1,0){14}}
\end{picture} 
& 217 & 5* & 11 & 1\\

13 &
\begin{picture}(140,40)
  \put(3,0){\circle*{6}}
  \put(23,0){\circle*{6}}
  \put(43,0){\circle{6}}
  \put(63,0){\circle*{6}}
  \put(83,0){\circle*{6}}
  \put(63,20){\circle{6}}
  \put(6,0){\line(1,0){14}}
  \put(26,0){\line(1,0){14}}
  \put(46,0){\line(1,0){14}}
  \put(66,0){\line(1,0){14}}
  \put(63,3){\line(0,1){14}}
\put(103,0){\circle*{6}}
\put(86,0){\line(1,0){14}}
\put(123,0){\circle*{6}}
\put(106,0){\line(1,0){14}}
\end{picture} 
& 111 & 1 & 10 & 0\\

14 &
\begin{picture}(140,40)
  \put(3,0){\circle*{6}}
  \put(23,0){\circle*{6}}
  \put(43,0){\circle*{6}}
  \put(63,0){\circle{6}}
  \put(83,0){\circle*{6}}
  \put(63,20){\circle*{6}}
  \put(6,0){\line(1,0){14}}
  \put(26,0){\line(1,0){14}}
  \put(46,0){\line(1,0){14}}
  \put(66,0){\line(1,0){14}}
  \put(63,3){\line(0,1){14}}
\put(103,0){\circle*{6}}
\put(86,0){\line(1,0){14}}
\put(123,0){\circle*{6}}
\put(106,0){\line(1,0){14}}
\end{picture} 
& 238 & 7* & 9 & 1\\

15 &
\begin{picture}(140,40)
  \put(3,0){\circle*{6}}
  \put(23,0){\circle*{6}}
  \put(43,0){\circle{6}}
  \put(63,0){\circle*{6}}
  \put(83,0){\circle*{6}}
  \put(63,20){\circle*{6}}
  \put(6,0){\line(1,0){14}}
  \put(26,0){\line(1,0){14}}
  \put(46,0){\line(1,0){14}}
  \put(66,0){\line(1,0){14}}
  \put(63,3){\line(0,1){14}}
\put(103,0){\circle*{6}}
\put(86,0){\line(1,0){14}}
\put(123,0){\circle*{6}}
\put(106,0){\line(1,0){14}}
\end{picture} 
& 159 & 1 & 9 & 0\\

16 &
\begin{picture}(140,40)
  \put(3,0){\circle*{6}}
  \put(23,0){\circle{6}}
  \put(43,0){\circle*{6}}
  \put(63,0){\circle*{6}}
  \put(83,0){\circle*{6}}
  \put(63,20){\circle*{6}}
  \put(6,0){\line(1,0){14}}
  \put(26,0){\line(1,0){14}}
  \put(46,0){\line(1,0){14}}
  \put(66,0){\line(1,0){14}}
  \put(63,3){\line(0,1){14}}
\put(103,0){\circle*{6}}
\put(86,0){\line(1,0){14}}
\put(123,0){\circle*{6}}
\put(106,0){\line(1,0){14}}
\end{picture} 
& 159 & 1 & 9 & 0\\

17 &
\begin{picture}(140,40)
  \put(3,0){\circle*{6}}
  \put(23,0){\circle*{6}}
  \put(43,0){\circle*{6}}
  \put(63,0){\circle*{6}}
  \put(83,0){\circle*{6}}
  \put(63,20){\circle*{6}}
  \put(6,0){\line(1,0){14}}
  \put(26,0){\line(1,0){14}}
  \put(46,0){\line(1,0){14}}
  \put(66,0){\line(1,0){14}}
  \put(63,3){\line(0,1){14}}
\put(103,0){\circle*{6}}
\put(86,0){\line(1,0){14}}
\put(123,0){\circle{6}}
\put(106,0){\line(1,0){14}}
\end{picture} 
& 127 & 1 & 8 & 0\\

\hline
\end{longtable}

\begin{longtable}{|c|c|c|c|c|c|}
\caption{N-regular inner automorphisms of $E_8$.}\label{tab:prinE8}
\endfirsthead
\hline
\multicolumn{6}{|l|}{\small\slshape N-regular automorphisms of $E_8$.} \\
\hline 
\endhead
\hline 
\endfoot
\endlastfoot

\hline

order & Kac diagram & \# orbits & \# components & $\dim$ & $\rank$ \\
\hline

2 &
\begin{picture}(150,40)
  \put(3,0){\circle*{6}}
  \put(23,0){\circle{6}}
  \put(43,0){\circle{6}}
  \put(63,0){\circle{6}}
  \put(83,0){\circle{6}}
  \put(43,20){\circle{6}}
  \put(6,0){\line(1,0){14}}
  \put(26,0){\line(1,0){14}}
  \put(46,0){\line(1,0){14}}
  \put(66,0){\line(1,0){14}}
  \put(43,3){\line(0,1){14}}
\put(103,0){\circle{6}}
\put(86,0){\line(1,0){14}}
\put(123,0){\circle{6}}
\put(106,0){\line(1,0){14}}
\put(143,0){\circle{6}}
\put(126,0){\line(1,0){14}}
\end{picture} 
& 115 & 1 & 120 & 8\\

3 &
\begin{picture}(150,40)
  \put(3,0){\circle{6}}
  \put(23,0){\circle{6}}
  \put(43,0){\circle{6}}
  \put(63,0){\circle{6}}
  \put(83,0){\circle{6}}
  \put(43,20){\circle*{6}}
  \put(6,0){\line(1,0){14}}
  \put(26,0){\line(1,0){14}}
  \put(46,0){\line(1,0){14}}
  \put(66,0){\line(1,0){14}}
  \put(43,3){\line(0,1){14}}
\put(103,0){\circle{6}}
\put(86,0){\line(1,0){14}}
\put(123,0){\circle{6}}
\put(106,0){\line(1,0){14}}
\put(143,0){\circle{6}}
\put(126,0){\line(1,0){14}}
\end{picture} 
& 101 & 1 & 80 & 4\\

4 &
\begin{picture}(150,40)
  \put(3,0){\circle{6}}
  \put(23,0){\circle{6}}
  \put(43,0){\circle{6}}
  \put(63,0){\circle{6}}
  \put(83,0){\circle*{6}}
  \put(43,20){\circle{6}}
  \put(6,0){\line(1,0){14}}
  \put(26,0){\line(1,0){14}}
  \put(46,0){\line(1,0){14}}
  \put(66,0){\line(1,0){14}}
  \put(43,3){\line(0,1){14}}
\put(103,0){\circle{6}}
\put(86,0){\line(1,0){14}}
\put(123,0){\circle{6}}
\put(106,0){\line(1,0){14}}
\put(143,0){\circle{6}}
\put(126,0){\line(1,0){14}}
\end{picture} 
& 144 & 2* & 60 & 4\\

5 &
\begin{picture}(150,40)
  \put(3,0){\circle{6}}
  \put(23,0){\circle{6}}
  \put(43,0){\circle{6}}
  \put(63,0){\circle*{6}}
  \put(83,0){\circle{6}}
  \put(43,20){\circle{6}}
  \put(6,0){\line(1,0){14}}
  \put(26,0){\line(1,0){14}}
  \put(46,0){\line(1,0){14}}
  \put(66,0){\line(1,0){14}}
  \put(43,3){\line(0,1){14}}
\put(103,0){\circle{6}}
\put(86,0){\line(1,0){14}}
\put(123,0){\circle{6}}
\put(106,0){\line(1,0){14}}
\put(143,0){\circle{6}}
\put(126,0){\line(1,0){14}}
\end{picture} 
& 105 & 1 & 48 & 2\\

6 &
\begin{picture}(150,40)
  \put(3,0){\circle{6}}
  \put(23,0){\circle{6}}
  \put(43,0){\circle{6}}
  \put(63,0){\circle*{6}}
  \put(83,0){\circle{6}}
  \put(43,20){\circle{6}}
  \put(6,0){\line(1,0){14}}
  \put(26,0){\line(1,0){14}}
  \put(46,0){\line(1,0){14}}
  \put(66,0){\line(1,0){14}}
  \put(43,3){\line(0,1){14}}
\put(103,0){\circle{6}}
\put(86,0){\line(1,0){14}}
\put(123,0){\circle{6}}
\put(106,0){\line(1,0){14}}
\put(143,0){\circle*{6}}
\put(126,0){\line(1,0){14}}
\end{picture} 
& 270 & 7* & 40 & 4\\

7 &
\begin{picture}(150,40)
  \put(3,0){\circle{6}}
  \put(23,0){\circle{6}}
  \put(43,0){\circle{6}}
  \put(63,0){\circle*{6}}
  \put(83,0){\circle{6}}
  \put(43,20){\circle{6}}
  \put(6,0){\line(1,0){14}}
  \put(26,0){\line(1,0){14}}
  \put(46,0){\line(1,0){14}}
  \put(66,0){\line(1,0){14}}
  \put(43,3){\line(0,1){14}}
\put(103,0){\circle{6}}
\put(86,0){\line(1,0){14}}
\put(123,0){\circle*{6}}
\put(106,0){\line(1,0){14}}
\put(143,0){\circle{6}}
\put(126,0){\line(1,0){14}}
\end{picture} 
& 144 & 1 & 35 & 1\\

8 &
\begin{picture}(150,40)
  \put(3,0){\circle{6}}
  \put(23,0){\circle{6}}
  \put(43,0){\circle*{6}}
  \put(63,0){\circle{6}}
  \put(83,0){\circle{6}}
  \put(43,20){\circle{6}}
  \put(6,0){\line(1,0){14}}
  \put(26,0){\line(1,0){14}}
  \put(46,0){\line(1,0){14}}
  \put(66,0){\line(1,0){14}}
  \put(43,3){\line(0,1){14}}
\put(103,0){\circle{6}}
\put(86,0){\line(1,0){14}}
\put(123,0){\circle*{6}}
\put(106,0){\line(1,0){14}}
\put(143,0){\circle{6}}
\put(126,0){\line(1,0){14}}
\end{picture} 
& 219 & 2* & 30 & 2\\

9 &
\begin{picture}(150,40)
  \put(3,0){\circle{6}}
  \put(23,0){\circle{6}}
  \put(43,0){\circle*{6}}
  \put(63,0){\circle{6}}
  \put(83,0){\circle{6}}
  \put(43,20){\circle{6}}
  \put(6,0){\line(1,0){14}}
  \put(26,0){\line(1,0){14}}
  \put(46,0){\line(1,0){14}}
  \put(66,0){\line(1,0){14}}
  \put(43,3){\line(0,1){14}}
\put(103,0){\circle{6}}
\put(86,0){\line(1,0){14}}
\put(123,0){\circle*{6}}
\put(106,0){\line(1,0){14}}
\put(143,0){\circle*{6}}
\put(126,0){\line(1,0){14}}
\end{picture} 
& 206 & 2* & 28 & 1\\

10 &
\begin{picture}(150,40)
  \put(3,0){\circle{6}}
  \put(23,0){\circle{6}}
  \put(43,0){\circle*{6}}
  \put(63,0){\circle{6}}
  \put(83,0){\circle{6}}
  \put(43,20){\circle{6}}
  \put(6,0){\line(1,0){14}}
  \put(26,0){\line(1,0){14}}
  \put(46,0){\line(1,0){14}}
  \put(66,0){\line(1,0){14}}
  \put(43,3){\line(0,1){14}}
\put(103,0){\circle*{6}}
\put(86,0){\line(1,0){14}}
\put(123,0){\circle{6}}
\put(106,0){\line(1,0){14}}
\put(143,0){\circle*{6}}
\put(126,0){\line(1,0){14}}
\end{picture} 
& 300 & 7* & 24 & 2\\

11 &
\begin{picture}(150,40)
  \put(3,0){\circle*{6}}
  \put(23,0){\circle{6}}
  \put(43,0){\circle*{6}}
  \put(63,0){\circle{6}}
  \put(83,0){\circle{6}}
  \put(43,20){\circle{6}}
  \put(6,0){\line(1,0){14}}
  \put(26,0){\line(1,0){14}}
  \put(46,0){\line(1,0){14}}
  \put(66,0){\line(1,0){14}}
  \put(43,3){\line(0,1){14}}
\put(103,0){\circle*{6}}
\put(86,0){\line(1,0){14}}
\put(123,0){\circle{6}}
\put(106,0){\line(1,0){14}}
\put(143,0){\circle{6}}
\put(126,0){\line(1,0){14}}
\end{picture} 
& 167 & 1 & 23 & 0\\

12 &
\begin{picture}(150,40)
  \put(3,0){\circle*{6}}
  \put(23,0){\circle{6}}
  \put(43,0){\circle*{6}}
  \put(63,0){\circle{6}}
  \put(83,0){\circle{6}}
  \put(43,20){\circle{6}}
  \put(6,0){\line(1,0){14}}
  \put(26,0){\line(1,0){14}}
  \put(46,0){\line(1,0){14}}
  \put(66,0){\line(1,0){14}}
  \put(43,3){\line(0,1){14}}
\put(103,0){\circle*{6}}
\put(86,0){\line(1,0){14}}
\put(123,0){\circle{6}}
\put(106,0){\line(1,0){14}}
\put(143,0){\circle*{6}}
\put(126,0){\line(1,0){14}}
\end{picture} 
& 398 & 10* & 20 & 2\\

13 &
\begin{picture}(150,40)
  \put(3,0){\circle{6}}
  \put(23,0){\circle{6}}
  \put(43,0){\circle*{6}}
  \put(63,0){\circle{6}}
  \put(83,0){\circle*{6}}
  \put(43,20){\circle{6}}
  \put(6,0){\line(1,0){14}}
  \put(26,0){\line(1,0){14}}
  \put(46,0){\line(1,0){14}}
  \put(66,0){\line(1,0){14}}
  \put(43,3){\line(0,1){14}}
\put(103,0){\circle{6}}
\put(86,0){\line(1,0){14}}
\put(123,0){\circle*{6}}
\put(106,0){\line(1,0){14}}
\put(143,0){\circle*{6}}
\put(126,0){\line(1,0){14}}
\end{picture} 
& 165 & 1 & 19 & 0\\

14 &
\begin{picture}(150,40)
  \put(3,0){\circle*{6}}
  \put(23,0){\circle{6}}
  \put(43,0){\circle*{6}}
  \put(63,0){\circle{6}}
  \put(83,0){\circle{6}}
  \put(43,20){\circle{6}}
  \put(6,0){\line(1,0){14}}
  \put(26,0){\line(1,0){14}}
  \put(46,0){\line(1,0){14}}
  \put(66,0){\line(1,0){14}}
  \put(43,3){\line(0,1){14}}
\put(103,0){\circle*{6}}
\put(86,0){\line(1,0){14}}
\put(123,0){\circle*{6}}
\put(106,0){\line(1,0){14}}
\put(143,0){\circle*{6}}
\put(126,0){\line(1,0){14}}
\end{picture} 
& 333 & 4* & 18 & 1\\

15 &
\begin{picture}(150,40)
  \put(3,0){\circle*{6}}
  \put(23,0){\circle{6}}
  \put(43,0){\circle*{6}}
  \put(63,0){\circle{6}}
  \put(83,0){\circle*{6}}
  \put(43,20){\circle{6}}
  \put(6,0){\line(1,0){14}}
  \put(26,0){\line(1,0){14}}
  \put(46,0){\line(1,0){14}}
  \put(66,0){\line(1,0){14}}
  \put(43,3){\line(0,1){14}}
\put(103,0){\circle{6}}
\put(86,0){\line(1,0){14}}
\put(123,0){\circle*{6}}
\put(106,0){\line(1,0){14}}
\put(143,0){\circle*{6}}
\put(126,0){\line(1,0){14}}
\end{picture} 
& 354 & 5* & 16 & 1\\

16 &
\begin{picture}(150,40)
  \put(3,0){\circle{6}}
  \put(23,0){\circle*{6}}
  \put(43,0){\circle{6}}
  \put(63,0){\circle*{6}}
  \put(83,0){\circle{6}}
  \put(43,20){\circle*{6}}
  \put(6,0){\line(1,0){14}}
  \put(26,0){\line(1,0){14}}
  \put(46,0){\line(1,0){14}}
  \put(66,0){\line(1,0){14}}
  \put(43,3){\line(0,1){14}}
\put(103,0){\circle*{6}}
\put(86,0){\line(1,0){14}}
\put(123,0){\circle{6}}
\put(106,0){\line(1,0){14}}
\put(143,0){\circle*{6}}
\put(126,0){\line(1,0){14}}
\end{picture} 
& 264 & 1 & 16 & 0\\

17 &
\begin{picture}(150,40)
  \put(3,0){\circle*{6}}
  \put(23,0){\circle{6}}
  \put(43,0){\circle*{6}}
  \put(63,0){\circle{6}}
  \put(83,0){\circle*{6}}
  \put(43,20){\circle*{6}}
  \put(6,0){\line(1,0){14}}
  \put(26,0){\line(1,0){14}}
  \put(46,0){\line(1,0){14}}
  \put(66,0){\line(1,0){14}}
  \put(43,3){\line(0,1){14}}
\put(103,0){\circle{6}}
\put(86,0){\line(1,0){14}}
\put(123,0){\circle*{6}}
\put(106,0){\line(1,0){14}}
\put(143,0){\circle{6}}
\put(126,0){\line(1,0){14}}
\end{picture} 
& 179 & 1 & 15 & 0\\

18 &
\begin{picture}(150,40)
  \put(3,0){\circle*{6}}
  \put(23,0){\circle*{6}}
  \put(43,0){\circle{6}}
  \put(63,0){\circle*{6}}
  \put(83,0){\circle{6}}
  \put(43,20){\circle*{6}}
  \put(6,0){\line(1,0){14}}
  \put(26,0){\line(1,0){14}}
  \put(46,0){\line(1,0){14}}
  \put(66,0){\line(1,0){14}}
  \put(43,3){\line(0,1){14}}
\put(103,0){\circle*{6}}
\put(86,0){\line(1,0){14}}
\put(123,0){\circle{6}}
\put(106,0){\line(1,0){14}}
\put(143,0){\circle*{6}}
\put(126,0){\line(1,0){14}}
\end{picture} 
& 397 & 5* & 14 & 1\\

19 &
\begin{picture}(150,40)
  \put(3,0){\circle*{6}}
  \put(23,0){\circle*{6}}
  \put(43,0){\circle*{6}}
  \put(63,0){\circle{6}}
  \put(83,0){\circle*{6}}
  \put(43,20){\circle{6}}
  \put(6,0){\line(1,0){14}}
  \put(26,0){\line(1,0){14}}
  \put(46,0){\line(1,0){14}}
  \put(66,0){\line(1,0){14}}
  \put(43,3){\line(0,1){14}}
\put(103,0){\circle{6}}
\put(86,0){\line(1,0){14}}
\put(123,0){\circle*{6}}
\put(106,0){\line(1,0){14}}
\put(143,0){\circle*{6}}
\put(126,0){\line(1,0){14}}
\end{picture} 
& 199 & 1 & 13 & 0\\

20 &
\begin{picture}(150,40)
  \put(3,0){\circle*{6}}
  \put(23,0){\circle*{6}}
  \put(43,0){\circle{6}}
  \put(63,0){\circle*{6}}
  \put(83,0){\circle{6}}
  \put(43,20){\circle*{6}}
  \put(6,0){\line(1,0){14}}
  \put(26,0){\line(1,0){14}}
  \put(46,0){\line(1,0){14}}
  \put(66,0){\line(1,0){14}}
  \put(43,3){\line(0,1){14}}
\put(103,0){\circle*{6}}
\put(86,0){\line(1,0){14}}
\put(123,0){\circle*{6}}
\put(106,0){\line(1,0){14}}
\put(143,0){\circle*{6}}
\put(126,0){\line(1,0){14}}
\end{picture} 
& 438 & 7* & 12 & 1\\

21 &
\begin{picture}(150,40)
  \put(3,0){\circle*{6}}
  \put(23,0){\circle{6}}
  \put(43,0){\circle*{6}}
  \put(63,0){\circle{6}}
  \put(83,0){\circle*{6}}
  \put(43,20){\circle*{6}}
  \put(6,0){\line(1,0){14}}
  \put(26,0){\line(1,0){14}}
  \put(46,0){\line(1,0){14}}
  \put(66,0){\line(1,0){14}}
  \put(43,3){\line(0,1){14}}
\put(103,0){\circle*{6}}
\put(86,0){\line(1,0){14}}
\put(123,0){\circle*{6}}
\put(106,0){\line(1,0){14}}
\put(143,0){\circle*{6}}
\put(126,0){\line(1,0){14}}
\end{picture} 
& 287 & 1 & 12 & 0\\

22 &
\begin{picture}(150,40)
  \put(3,0){\circle{6}}
  \put(23,0){\circle*{6}}
  \put(43,0){\circle{6}}
  \put(63,0){\circle*{6}}
  \put(83,0){\circle*{6}}
  \put(43,20){\circle*{6}}
  \put(6,0){\line(1,0){14}}
  \put(26,0){\line(1,0){14}}
  \put(46,0){\line(1,0){14}}
  \put(66,0){\line(1,0){14}}
  \put(43,3){\line(0,1){14}}
\put(103,0){\circle*{6}}
\put(86,0){\line(1,0){14}}
\put(123,0){\circle*{6}}
\put(106,0){\line(1,0){14}}
\put(143,0){\circle*{6}}
\put(126,0){\line(1,0){14}}
\end{picture} 
& 319 & 1 & 12 & 0\\

23 &
\begin{picture}(150,40)
  \put(3,0){\circle*{6}}
  \put(23,0){\circle{6}}
  \put(43,0){\circle*{6}}
  \put(63,0){\circle*{6}}
  \put(83,0){\circle*{6}}
  \put(43,20){\circle{6}}
  \put(6,0){\line(1,0){14}}
  \put(26,0){\line(1,0){14}}
  \put(46,0){\line(1,0){14}}
  \put(66,0){\line(1,0){14}}
  \put(43,3){\line(0,1){14}}
\put(103,0){\circle*{6}}
\put(86,0){\line(1,0){14}}
\put(123,0){\circle*{6}}
\put(106,0){\line(1,0){14}}
\put(143,0){\circle*{6}}
\put(126,0){\line(1,0){14}}
\end{picture} 
& 233 & 1 & 11 & 0\\

24 &
\begin{picture}(150,40)
  \put(3,0){\circle*{6}}
  \put(23,0){\circle*{6}}
  \put(43,0){\circle{6}}
  \put(63,0){\circle*{6}}
  \put(83,0){\circle*{6}}
  \put(43,20){\circle*{6}}
  \put(6,0){\line(1,0){14}}
  \put(26,0){\line(1,0){14}}
  \put(46,0){\line(1,0){14}}
  \put(66,0){\line(1,0){14}}
  \put(43,3){\line(0,1){14}}
\put(103,0){\circle*{6}}
\put(86,0){\line(1,0){14}}
\put(123,0){\circle*{6}}
\put(106,0){\line(1,0){14}}
\put(143,0){\circle*{6}}
\put(126,0){\line(1,0){14}}
\end{picture} 
& 478 & 8* & 10 & 1\\

25 &
\begin{picture}(150,40)
  \put(3,0){\circle*{6}}
  \put(23,0){\circle*{6}}
  \put(43,0){\circle*{6}}
  \put(63,0){\circle{6}}
  \put(83,0){\circle*{6}}
  \put(43,20){\circle*{6}}
  \put(6,0){\line(1,0){14}}
  \put(26,0){\line(1,0){14}}
  \put(46,0){\line(1,0){14}}
  \put(66,0){\line(1,0){14}}
  \put(43,3){\line(0,1){14}}
\put(103,0){\circle*{6}}
\put(86,0){\line(1,0){14}}
\put(123,0){\circle*{6}}
\put(106,0){\line(1,0){14}}
\put(143,0){\circle*{6}}
\put(126,0){\line(1,0){14}}
\end{picture} 
& 319 & 1 & 10 & 0\\

26 &
\begin{picture}(150,40)
  \put(3,0){\circle*{6}}
  \put(23,0){\circle*{6}}
  \put(43,0){\circle*{6}}
  \put(63,0){\circle*{6}}
  \put(83,0){\circle{6}}
  \put(43,20){\circle*{6}}
  \put(6,0){\line(1,0){14}}
  \put(26,0){\line(1,0){14}}
  \put(46,0){\line(1,0){14}}
  \put(66,0){\line(1,0){14}}
  \put(43,3){\line(0,1){14}}
\put(103,0){\circle*{6}}
\put(86,0){\line(1,0){14}}
\put(123,0){\circle*{6}}
\put(106,0){\line(1,0){14}}
\put(143,0){\circle*{6}}
\put(126,0){\line(1,0){14}}
\end{picture} 
& 319 & 1 & 10 & 0\\

27 &
\begin{picture}(150,40)
  \put(3,0){\circle*{6}}
  \put(23,0){\circle*{6}}
  \put(43,0){\circle*{6}}
  \put(63,0){\circle*{6}}
  \put(83,0){\circle*{6}}
  \put(43,20){\circle*{6}}
  \put(6,0){\line(1,0){14}}
  \put(26,0){\line(1,0){14}}
  \put(46,0){\line(1,0){14}}
  \put(66,0){\line(1,0){14}}
  \put(43,3){\line(0,1){14}}
\put(103,0){\circle{6}}
\put(86,0){\line(1,0){14}}
\put(123,0){\circle*{6}}
\put(106,0){\line(1,0){14}}
\put(143,0){\circle*{6}}
\put(126,0){\line(1,0){14}}
\end{picture} 
& 319 & 1 & 10 & 0\\

28 &
\begin{picture}(150,40)
  \put(3,0){\circle*{6}}
  \put(23,0){\circle*{6}}
  \put(43,0){\circle*{6}}
  \put(63,0){\circle*{6}}
  \put(83,0){\circle*{6}}
  \put(43,20){\circle*{6}}
  \put(6,0){\line(1,0){14}}
  \put(26,0){\line(1,0){14}}
  \put(46,0){\line(1,0){14}}
  \put(66,0){\line(1,0){14}}
  \put(43,3){\line(0,1){14}}
\put(103,0){\circle*{6}}
\put(86,0){\line(1,0){14}}
\put(123,0){\circle{6}}
\put(106,0){\line(1,0){14}}
\put(143,0){\circle*{6}}
\put(126,0){\line(1,0){14}}
\end{picture} 
& 319 & 1 & 10 & 0\\

29 &
\begin{picture}(150,40)
  \put(3,0){\circle*{6}}
  \put(23,0){\circle*{6}}
  \put(43,0){\circle*{6}}
  \put(63,0){\circle*{6}}
  \put(83,0){\circle*{6}}
  \put(43,20){\circle*{6}}
  \put(6,0){\line(1,0){14}}
  \put(26,0){\line(1,0){14}}
  \put(46,0){\line(1,0){14}}
  \put(66,0){\line(1,0){14}}
  \put(43,3){\line(0,1){14}}
\put(103,0){\circle*{6}}
\put(86,0){\line(1,0){14}}
\put(123,0){\circle*{6}}
\put(106,0){\line(1,0){14}}
\put(143,0){\circle{6}}
\put(126,0){\line(1,0){14}}
\end{picture} 
& 255 & 1 & 9 & 0\\

\hline
\end{longtable}

\begin{longtable}{|c|c|c|c|c|c|}
\caption{N-regular inner automorphisms of $F_4$.}\label{tab:prinF4}
\endfirsthead
\hline
\multicolumn{6}{|l|}{\small\slshape N-regular automorphisms of $F_4$.} \\
\hline 
\endhead
\hline 
\endfoot
\endlastfoot

\hline

order & Kac diagram & \# orbits & \# components & $\dim$ & $\rank$ \\
\hline

2 & 
\begin{picture}(130,15)
  \put(10,5){\circle{6}}
  \put(40,5){\circle*{6}}
  \put(70,5){\circle{6}}
  \put(100,5){\circle{6}}
  \put(130,5){\circle{6}}
  \put(13,5){\line(1,0){24}}
  \put(43,5){\line(1,0){24}}
  \put(72,7){\line(1,0){26}}
  \put(72,3){\line(1,0){26}}
  \put(80,1){\Large $>$}
  \put(103,5){\line(1,0){24}}
\end{picture}
& 26 & 1 & 24 & 4\\

3 & 
\begin{picture}(130,15)
  \put(10,5){\circle{6}}
  \put(40,5){\circle{6}}
  \put(70,5){\circle*{6}}
  \put(100,5){\circle{6}}
  \put(130,5){\circle{6}}
  \put(13,5){\line(1,0){24}}
  \put(43,5){\line(1,0){24}}
  \put(72,7){\line(1,0){26}}
  \put(72,3){\line(1,0){26}}
  \put(80,1){\Large $>$}
  \put(103,5){\line(1,0){24}}
\end{picture}
& 19 & 1 & 16 & 2
\\

4 & 
\begin{picture}(130,15)
  \put(10,5){\circle*{6}}
  \put(40,5){\circle{6}}
  \put(70,5){\circle*{6}}
  \put(100,5){\circle{6}}
  \put(130,5){\circle{6}}
  \put(13,5){\line(1,0){24}}
  \put(43,5){\line(1,0){24}}
  \put(72,7){\line(1,0){26}}
  \put(72,3){\line(1,0){26}}
  \put(80,1){\Large $>$}
  \put(103,5){\line(1,0){24}}
\end{picture}
& 29 & 3* & 12 & 2 \\

5 & 
\begin{picture}(130,15)
  \put(10,5){\circle{6}}
  \put(40,5){\circle{6}}
  \put(70,5){\circle*{6}}
  \put(100,5){\circle{6}}
  \put(130,5){\circle*{6}}
  \put(13,5){\line(1,0){24}}
  \put(43,5){\line(1,0){24}}
  \put(72,7){\line(1,0){26}}
  \put(72,3){\line(1,0){26}}
  \put(80,1){\Large $>$}
  \put(103,5){\line(1,0){24}}
\end{picture}
& 15 & 1 & 11 & 0\\

6 & 
\begin{picture}(130,15)
  \put(10,5){\circle*{6}}
  \put(40,5){\circle{6}}
  \put(70,5){\circle*{6}}
  \put(100,5){\circle{6}}
  \put(130,5){\circle*{6}}
  \put(13,5){\line(1,0){24}}
  \put(43,5){\line(1,0){24}}
  \put(72,7){\line(1,0){26}}
  \put(72,3){\line(1,0){26}}
  \put(80,1){\Large $>$}
  \put(103,5){\line(1,0){24}}
\end{picture}
& 35 & 6* & 8 & 2\\

7 & 
\begin{picture}(130,15)
  \put(10,5){\circle*{6}}
  \put(40,5){\circle*{6}}
  \put(70,5){\circle{6}}
  \put(100,5){\circle*{6}}
  \put(130,5){\circle{6}}
  \put(13,5){\line(1,0){24}}
  \put(43,5){\line(1,0){24}}
  \put(72,7){\line(1,0){26}}
  \put(72,3){\line(1,0){26}}
  \put(80,1){\Large $>$}
  \put(103,5){\line(1,0){24}}
\end{picture}
& 13 & 1 & 7 & 0\\

8 & 
\begin{picture}(130,15)
  \put(10,5){\circle*{6}}
  \put(40,5){\circle*{6}}
  \put(70,5){\circle*{6}}
  \put(100,5){\circle{6}}
  \put(130,5){\circle*{6}}
  \put(13,5){\line(1,0){24}}
  \put(43,5){\line(1,0){24}}
  \put(72,7){\line(1,0){26}}
  \put(72,3){\line(1,0){26}}
  \put(80,1){\Large $>$}
  \put(103,5){\line(1,0){24}}
\end{picture}
& 30 & 4* & 6 & 1\\

9 & 
\begin{picture}(130,15)
  \put(10,5){\circle*{6}}
  \put(40,5){\circle*{6}}
  \put(70,5){\circle{6}}
  \put(100,5){\circle*{6}}
  \put(130,5){\circle*{6}}
  \put(13,5){\line(1,0){24}}
  \put(43,5){\line(1,0){24}}
  \put(72,7){\line(1,0){26}}
  \put(72,3){\line(1,0){26}}
  \put(80,1){\Large $>$}
  \put(103,5){\line(1,0){24}}
\end{picture}
& 19 & 1 & 6 & 0\\

10 & 
\begin{picture}(130,15)
  \put(10,5){\circle*{6}}
  \put(40,5){\circle{6}}
  \put(70,5){\circle*{6}}
  \put(100,5){\circle*{6}}
  \put(130,5){\circle*{6}}
  \put(13,5){\line(1,0){24}}
  \put(43,5){\line(1,0){24}}
  \put(72,7){\line(1,0){26}}
  \put(72,3){\line(1,0){26}}
  \put(80,1){\Large $>$}
  \put(103,5){\line(1,0){24}}
\end{picture}
& 19 & 1 & 6 & 0\\

11 & 
\begin{picture}(130,15)
  \put(10,5){\circle{6}}
  \put(40,5){\circle*{6}}
  \put(70,5){\circle*{6}}
  \put(100,5){\circle*{6}}
  \put(130,5){\circle*{6}}
  \put(13,5){\line(1,0){24}}
  \put(43,5){\line(1,0){24}}
  \put(72,7){\line(1,0){26}}
  \put(72,3){\line(1,0){26}}
  \put(80,1){\Large $>$}
  \put(103,5){\line(1,0){24}}
\end{picture}
& 15 & 1 & 5 & 0
\\

\hline
\end{longtable}

\begin{longtable}{|c|c|c|c|c|c|}
\caption{N-regular inner automorphisms of $G_2$.}\label{tab:prinG2}
\endfirsthead
\hline
\multicolumn{6}{|l|}{\small\slshape N-regular automorphisms of $G_2$.} \\
\hline 
\endhead
\hline 
\endfoot
\endlastfoot

\hline

order & Kac diagram & \# orbits & \# components & $\dim$ & $\rank$ \\
\hline

2 & 
\begin{picture}(70,7)
  \put(5,0){\circle{6}}
  \put(35,0){\circle*{6}}
  \put(65,0){\circle{6}}
  \put(8,0){\line(1,0){24}}
  \put(35,-3){\line(1,0){30}}
  \put(38,0){\line(1,0){24}}
  \put(35,3){\line(1,0){30}}
  \put(45,-4){\Large $>$}
\end{picture}
& 5 & 1 & 6 & 2\\

3 & 
\begin{picture}(70,7)
  \put(5,0){\circle*{6}}
  \put(35,0){\circle*{6}}
  \put(65,0){\circle{6}}
  \put(8,0){\line(1,0){24}}
  \put(35,-3){\line(1,0){30}}
  \put(38,0){\line(1,0){24}}
  \put(35,3){\line(1,0){30}}
  \put(45,-4){\Large $>$}
\end{picture}
& 6 & 2* & 4 & 1\\

4 & 
\begin{picture}(70,7)
  \put(5,0){\circle*{6}}
  \put(35,0){\circle{6}}
  \put(65,0){\circle*{6}}
  \put(8,0){\line(1,0){24}}
  \put(35,-3){\line(1,0){30}}
  \put(38,0){\line(1,0){24}}
  \put(35,3){\line(1,0){30}}
  \put(45,-4){\Large $>$}
\end{picture}
& 4 & 1 & 4 & 0\\

5 & 
\begin{picture}(70,7)
  \put(5,0){\circle{6}}
  \put(35,0){\circle*{6}}
  \put(65,0){\circle*{6}}
  \put(8,0){\line(1,0){24}}
  \put(35,-3){\line(1,0){30}}
  \put(38,0){\line(1,0){24}}
  \put(35,3){\line(1,0){30}}
  \put(45,-4){\Large $>$}
\end{picture} 
& 3 & 1 & 3 & 0\\

\hline
\end{longtable}

\def\cprime{$'$} \def\cprime{$'$} \def\Dbar{\leavevmode\lower.6ex\hbox to
  0pt{\hskip-.23ex \accent"16\hss}D} \def\cprime{$'$} \def\cprime{$'$}
  \def\cprime{$'$}

\end{document}